\documentclass[11pt]{amsart} %
\usepackage{amssymb,amsfonts,amsmath}
\usepackage[dvips]{graphicx}

\usepackage{psfrag}
\newtheorem{theo}{Theorem}[section] 
\newtheorem{prop}[theo]{Proposition}
\newtheorem{lemma}[theo]{Lemma}
\newtheorem{coro}[theo]{Corollary}

\newcommand{\mc}[1]{\mathcal{#1}}

\newcommand{\nib}[1]{\noindent {\bf #1}}

\newcommand{\sm}{\setminus}

\newcommand{\eps}{\epsilon}

\newcommand{\sub}{\subseteq}

\def\qed{\hfill $\Box$}

\def\COMMENT#1{}
\def\TASK#1{}
\def\noproof{{\unskip\nobreak\hfill\penalty50\hskip2em\hbox{}\nobreak\hfill%
       $\square$\parfillskip=0pt\finalhyphendemerits=0\par}\goodbreak}
\def\endproof{\noproof\bigskip}
\newdimen\margin   
\def\textno#1&#2\par{%
   \margin=\hsize
   \advance\margin by -4\parindent
          \setbox1=\hbox{\sl#1}%
   \ifdim\wd1 < \margin
      $$\box1\eqno#2$$%
   \else
      \bigbreak
      \hbox to \hsize{\indent$\vcenter{\advance\hsize by -3\parindent
      \sl\noindent#1}\hfil#2$}%
      \bigbreak
   \fi}

\def\enddiscard{}
\long\def\discard#1\enddiscard{}

\title{Loose Hamilton Cycles in Hypergraphs}
\author{Peter Keevash \and Daniela K\"uhn \and Richard Mycroft \and Deryk Osthus}
\date{}
\thanks {D.~K\"uhn was partially supported by the EPSRC, grant no.~EP/D50564X/1.
D.~Osthus was partially supported by the EPSRC, grant no.~EP/E02162X/1. P.~Keevash was partially supported by the ERC, grant no.~239696, and by the EPSRC, grant no.~EP/G056730/1}

\begin{document}

\vspace*{-0.8cm}
\begin{abstract}
We prove that any $k$-uniform hypergraph on $n$ vertices with minimum degree at least $\frac{n}{2(k-1)}+o(n)$ 
contains a loose Hamilton cycle. The proof strategy is similar to that used by K\"uhn and Osthus for the 
3-uniform case. Though some additional difficulties arise in the $k$-uniform case, our argument here is 
considerably simplified by applying the recent hypergraph blow-up lemma of Keevash.
\end{abstract}

\maketitle

\vspace*{-0.6cm}
\section{Introduction}\label{intro}

A fundamental theorem of Dirac~\cite{DIRAC} states that any graph on $n$ vertices with minimum degree at least $n/2$ 
contains a Hamilton cycle. A natural question is whether this theorem can be extended to hypergraphs. 

For this, we first need to extend the notions of minimum degree and of Hamilton cycles to hypergraphs. 
A \emph{$k$-uniform hypergraph} or \emph{$k$-graph} $H$ consists of a vertex set $V$ and a set of edges  
each consisting of $k$ vertices. We will often identify $H$ with its edge set and write $e \in H$ if $e$ is an 
edge of~$H$. Given a $k$-graph $H$, we say that a set of $k-1$ vertices 
$T \in \binom{V}{k-1}$ has \emph{neighbourhood} $N_H(T) = \{x \in V: \{x\} \cup T \in H \}$. 
The \emph{degree of $T$} is $d_{k-1}(T) = |N_H(T)|$. The \emph{minimum degree} of $H$ is the minimum 
size of  such a neighbourhood, that is, $\delta_{k-1}(H) = \min \{d_{k-1}(T) : T \in \binom{V}{k-1}\}$. 

\begin{figure}\label{fig:cycles}
\centering\footnotesize
\includegraphics[width=0.4\columnwidth]{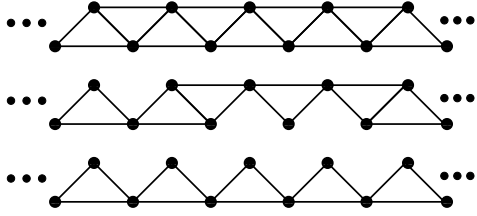}
\caption{Segments of a tight cycle (top), a generic  cycle (middle) and a loose cycle (bottom).}
\end{figure}

We say that a $k$-graph $C$ is a \emph{cycle of order $n$} if its vertices can be given a cyclic ordering 
$v_1,\dots, v_n$ so that every consecutive pair $v_i, v_{i+1}$ lies in an edge of $C$ and every edge 
of $C$ consists of $k$ consecutive vertices. A cycle of order $n$ is \emph{tight} if every set of $k$ 
consecutive vertices forms an edge; it is \emph{loose} if every pair of adjacent edges intersects in a 
single vertex, with the possible exception of one pair of edges, which may intersect in more than 
one vertex. This final condition allows us to consider loose cycles whose order is not a multiple 
of $k-1$. Figure~\ref{fig:cycles} shows the structure of each of these cycle types. 
A \emph{Hamilton cycle} 
in a $k$-graph $H$ is a sub-$k$-graph of $H$ which is a cycle containing every vertex of $H$. 

R\"odl, Ruci\'nski and Szemer\'edi~\cite{RRS,RRS2} showed that for any $\eta>0$ there is 
an $n_0$ so that if $n>n_0$ then any $k$-graph $H$ on $n$ vertices with 
minimum degree $\delta_{k-1}(H)\ge n/2+ \eta n$ contains a tight Hamilton cycle
(this improved an earlier bound by Katona and Kierstead~\cite{KK99}).
They gave a construction which shows that this result is best possible up to the error term $\eta n$.
In this paper, we prove the analogous result for loose Hamilton cycles.

\begin{theo} \label{main}
For all $k\geq 3$ and any $\eta>0$ there exists $n_0$ so that if $n>n_0$ then any $k$-graph $H$
on $n$ vertices with $\delta_{k-1}(H) > ( \frac{1}{2(k-1)} + \eta)n$ contains a loose Hamilton cycle.
\end{theo}

The case when $k=3$ was proved by K\"uhn and Osthus~\cite{KO}. 
We will use a similar method of proof for general $k$-graphs, but this will be greatly simplified 
by the use of the recent blow-up lemma of Keevash~\cite{K2}. 

Proposition~\ref{bestpos} shows that Theorem~\ref{main} is best possible up to the error term $\eta n$. 
In fact, Proposition~\ref{bestpos} actually tells us more than this, namely that up to the error term, 
this minimum degree condition is best possible to ensure the existence of any (not necessarily loose) 
Hamilton cycle in $H$. This means that the minimum degree needed to find a Hamilton cycle in a 
$k$-graph of order $n$ is $\frac{n}{2(k-1)} + o(n)$.

Whilst finalizing this paper we learnt that H\`an and Schacht~\cite{HS} independently and 
simultaneously proved Theorem~\ref{main}, using a different approach.
The result in~\cite{HS} also covers the notion of a $k$-uniform $\ell$-cycle
for $\ell<k/2$ (here one requires consecutive edges to intersect in precisely $\ell$ vertices). More recently K\"uhn, Mycroft and Osthus~\cite{KMO} further developed the method of H\`an and Schacht to include all $\ell$ such that $k-\ell \nmid k$ (the remaining values of $\ell$ are covered by the results of R\"odl, Ruci\'nski and Szemer\'edi~\cite{RRS, RRS2}).

There is also the notion of a \emph{Berge-cycle}, which consists of a sequence of vertices where each pair
of consecutive vertices is contained in a common edge. This is less restrictive than the cycles considered
in this paper. Hamiltonian Berge-cycles were studied in~\cite{bermond}.

\section{Extremal example and outline of the proof} \label{sec:extremal}

The next proposition shows that Theorem~\ref{main} is best possible,
up to the error term~$\eta n$. 

\begin{prop} \label{bestpos}
For all integers $k \geq 3$ and $n\ge 2k-1$, there exists a $k$-graph $H$ on $n$ vertices such that
$\delta_{k-1}(H) \geq \lceil\frac{n}{2k-2}\rceil-1$ but $H$ does not contain a Hamilton cycle.
\end{prop}

\nib{Proof.} Let $V_1$ and $V_2$ be disjoint sets of size $\lceil\frac{n}{2k-2}\rceil-1$ and
$n - \lceil\frac{n}{2k-2}\rceil+1$ respectively. Let $H$ be the $k$-graph on the vertex set
$V = V_1 \cup V_2$, with $e \in \binom{V}{k}$ an edge if and only if $e \cap V_1 \neq \emptyset$,
that is, if $e$ contains at least one vertex from $V_1$. Then $H$ has minimum degree
$\delta_{k-1}(H) = \lceil\frac{n}{2k-2}\rceil-1$. However, any cyclic ordering of the vertices of~$H$
must contain $2k-2$ consecutive vertices $v_1,\dots,v_{2k-2}$ from $V_2$, but then $v_{k-1}$ and
$v_k$ cannot be contained in a common edge consisting of $k$ consecutive vertices, and so $H$ cannot contain a
Hamilton cycle. \qed

\medskip

In our proof of Theorem~\ref{main} we construct the loose Hamilton cycle by finding several paths and 
joining them into a spanning cycle. 
Here a $k$-graph $P$ is a \emph{path} if its vertices can be given a linear ordering 
such that every edge of $P$ consists of $k$ consecutive vertices, and so that every pair of consecutive 
vertices of $P$ lie in an edge of $P$. Similarly as for cycles, we say that a path $P$ is \emph{loose} 
if edges of $P$ intersect in at most one vertex. The ordering of the vertices of $P$ naturally 
gives an ordering of the edges of $P$. We say that any vertex of $P$ which lies in the initial edge of 
$P$, but not the second edge of $P$, is an \emph{initial vertex}. Similarly, any vertex of $P$ which lies in the 
final edge of $P$ but not the penultimate edge is a \emph{final vertex}. Also, we refer to vertices of 
$P$ which lie in more than one edge of $P$ as \emph{link vertices}. Thus, for example, a loose path $P$ has 
$k-1$ initial vertices, $k-1$ final vertices, and one link vertex in each pair of consecutive edges.

In Section~\ref{sec:regularity}, we shall introduce various ideas we will need in the proof of Theorem~\ref{main}.
In particular, we will state a version of the hypergraph regularity lemma due to R\"odl and Schacht~\cite{RS}
and Theorem~\ref{robust-universal} due to Keevash~\cite{K2}.
The latter provides a useful way of applying the hypergraph blow-up lemma. In Section~\ref{sec:prelims}, we shall prove various 
auxiliary results, including a result on finding loose paths in complete $k$-partite $k$-graphs, 
and an approximate minimum degree condition to guarantee a near-perfect packing of $H$ with a particular 
$k$-graph $\mathcal{A}_k$. Finally, in Section~\ref{proof} we shall prove Theorem~\ref{main} as follows.%
     \COMMENT{changed pf later on a bit, so overview looks a bit different now}

\subsection{Imposing structure on $H$.} In Section~\ref{structure} we use the hypergraph regularity lemma
to split $H$ into 
$k$-partite $k$-graphs $H^i$ on disjoint vertex sets $X^i$. 
These $k$-graphs $H^i$ will be suitable for embedding almost spanning 
loose paths, and all the vertices of $H$ not contained in 
any of the $X^i$ will be included in an `exceptional' loose path $L_e$
(actually, if $|V(H)|$ is not divisible by $k-1$, then $L_e$ will contain two consecutive edges which 
intersect in more than one vertex). The requirement that $H^i$ contains an almost spanning loose path
means that the vertex classes of the $H^i$ must have suitable size.
We achieve this by first defining a suitable `reduced $k$-graph' $R$ of $H$. 
Then we cover almost all vertices of $R$ by copies of a suitable auxiliary $k$-graph ${\mathcal A}_k$.
For each copy of ${\mathcal A}_k$, the corresponding sub-$k$-graph of $H$ is then split into the same number
of disjoint $H^i$.

\subsection{The linking strategy.}
In Section~\ref{linking} we shall use the structure imposed on $H$ to find a 
Hamilton cycle in $H$ by the following process.
\begin{itemize}
\item[(a)] The $k$-graphs $H^i$ are connected by means of a walk $W = e_1, \dots, e_\ell$ in the `supplementary graph'.
This graph (which we will define in Section~\ref{supp}) has vertices $1,\dots,t'$ corresponding to the $k$-graphs $H^i$.
\item[(b)] Using Lemma~\ref{interpath}, each edge $e_j$ of $W$ is used to create a short `connecting' loose 
path $L_j$ in $H$ joining two different $H^i$s.
\item[(c)] $L_e$ and the paths $L_j$ are extended to `prepaths' (these can be thought 
of as a path minus an initial vertex and a final vertex) $L_e^* = I_0L_eF_0$ and $L_j^* = I_jL_jF_j$, where 
$I_0, F_0$ and all $I_j, F_j$ are sets of size $k-2$. These prepaths have the property that 
there are large sets $I'_j$ and $F'_j$ such that $L_j^*$ 
can be extended to a loose path by adding any vertex of $I'_j$ as an initial vertex and any vertex of 
$F'_j$ as a final vertex. Similarly there are large sets $I'_{\ell+1}$ and $F'_0$ so that
$L_e^*$ can be extended to a path by adding any vertex of $I'_{\ell+1}$ as an initial vertex 
and any vertex of $F'_0$ as a final vertex. 
$I'_{j+1}$ and $F'_j$ both lie in the same $H^i$ (for all $j=0,\dots,\ell$).
\item[(d)] For each $H^i$ and for all those pairs $I'_{j+1},F'_j$ which lie in $H^i$, we choose a loose path 
$L'_{j+1}$ inside $H^i$ from $F'_j$ to $I'_{j+1}$. 
For each $i$, we will use the hypergraph blow-up lemma (in the form of Theorem~\ref{robust-universal})
to ensure that together all those $L'_j$ which lie in~$H^i$ use all the remaining vertices of $H^i$.
\item[(e)] The loose Hamilton cycle is then the concatenation $L_e^*L_1'L_1^*\dots L'_\ell L_\ell^* L'_{\ell+1}$.
\end{itemize}

\subsection{Controlling divisibility.}

Note that the number of vertices of a loose path is $1$ modulo $k-1$. So in order to apply 
Theorem~\ref{robust-universal} to obtain spanning loose paths in a subgraph of $H^i$, 
we need this subgraph to satisfy this condition.
So we choose our paths sequentially to satisfy the following congruences modulo $k-1$.
\begin{itemize}
\item[(a)] $L_e$ is chosen with $|V(H)\sm V(L_e)| \equiv -1$. 
\item[(b)]  Let $X^i(j-1)$ be the subset of $X^i$ obtained by removing $V(L_1),\dots,V(L_{j-1})$.
(All the $X^i$ will be disjoint from~$V(L_e)$.)
Let $d_i$ be the number of times that $W$ visits $H^i$.
When choosing $L_j$, for every $X^i$ it traverses (except the final one) 
we arrange to intersect $X^i(j-1)$  in a set of size 
$\equiv t_i(j) \equiv |X^i(j-1)|+d_i$ (the size modulo $k-1$ of the intersection of $L_j$ with 
the final $X^i$ it traverses is then determined
by the sizes of the other intersections). 
The choice of $L_e$ in (a) ensures that after all $L_j$ have been picked, 
the remaining part $X^i(\ell)$ of $X^i$ has size $\equiv -d_i$.
\item[(c)] Each $L_j$ is extended to a prepath $L_j^*$ by adding~$I_j$ and~$F_j$.
Similarly, $L_e$ is extended into a prepath $L_e^*$ by adding~$I_0$ and~$F_0$.
Now the remaining part of $X^i$ has size $\equiv d_i$.
\item[(d)] It remains to select $d_i$ paths $L'_j$ within each $X^i$: each uses $\equiv 1$ vertices,
so the divisibility conditions are satisfied.
\end{itemize}

\section{Regularity and the Blow-up Lemma} \label{sec:regularity}
\subsection{Graphs and complexes}

We begin with some notation. By $[r]$ we denote the set of integers from 1 to $r$. For a set $A$, we use 
$\binom{A}{k}$ to denote the collection of subsets of $A$ of size $k$, and similarly $\binom{A}{\leq k}$ 
to denote the collection of non-empty%
    \COMMENT{have non-empty here since that's what we need in most places, for example in the
def of $\eps$-regular} 
subsets of $A$ of size at most $k$. We write $x = y \pm z$ 
to mean that $y-z \leq x \leq y+z$. We shall omit floors and ceilings throughout this paper whenever they 
do not affect the argument.

A \emph{hypergraph} $H$ consists of a vertex set $V(H)$ and an edge set, such that each edge 
$e$ of the hypergraph satisfies $e \subseteq V(H)$. So a $k$-graph as defined in
Section~\ref{intro} is a
hypergraph in which all the edges are of size $k$. We say that a hypergraph $H$ is a \emph{$k$-complex} 
if every edge has size at most $k$ and $H$ forms a simplicial complex, that is, if
$e_1 \in H$ and $e_2 \subseteq e_1$ then $e_2 \in H$. As for $k$-graphs we identify a hypergraph $H$ with the set of its edges. So $|H|$ is the number of edges in~$H$, and if $G$ and $H$ are hypergraphs then $G \sm H$ is formed by removing from $G$ any edge which also lies in $H$. If $H$ is a hypergraph with vertex set $V$ then for any $V' \subseteq V$ the \emph{restriction $H[V']$ of $H$ to $V'$} is defined to have vertex set $V'$ and all edges of $H$ which are contained in $V'$ as edges. Also, for any hypergraphs $G$ and $H$ we define $G - H$ to be the hypergraph $G[V(G) \sm V(H)]$.

We say that a hypergraph $H$ is \emph{$r$-partite} 
if its vertex set $X$ is divided into $r$ pairwise-disjoint parts $X_1, \dots, X_r$, in such a way that for any edge 
$e \in H$, $|e \cap X_i| \leq 1$ for each $i$. 
We call the $X_i$ the \emph{vertex classes} of $H$
and say that the partition $X_1, \dots, X_r$ of $X$ is \emph{equitable} if all the $X_i$ have the same size. 
We say that a set $A \subseteq X$ is \emph{$r$-partite} if 
$|A \cap X_i| \leq 1$ for each $i$. So every edge of an $r$-partite hypergraph is $r$-partite. In the 
same way we may also speak of $r$-partite $k$-graphs and $r$-partite $k$-complexes. 
Given a $k$-graph $H$, we define a $k$-complex 
$H^\leq = \{e_1 \colon e_1 \subseteq e_2$ and $e_2 \in H\}$ and 
a $(k-1)$-complex $H^< = \{e_1 \colon e_1 \subset e_2$ and $e_2 \in H\}$.
Conversely, for a $k$-complex $H$ we define the $k$-graph $H_=$ to be the `top level' of 
$H$, i.e.~$H_= = \{e \in H \colon |e| = k\}$.
(Here $V(H)=V(H^\le)=V(H^<)=V(H_=)$.) 

Given a $k$-graph $G$ and a set~$W$ of vertices of $G$, we denote by $G[W]$ the sub-$k$-graph of $G$ 
obtained by removing all vertices and edges not contained in $W$
(in this case, we say $G$ is \emph{restricted to $W$}).
For a $k$-graph $G$ and a sub-$k$-graph $H \subseteq G$ write $G - H$ for $G[V(G)\sm V(H)]$.

Let $X_1, \dots, X_r$ be pairwise-disjoint sets of vertices, and let $X = X_1 \cup \dots \cup X_r$. Given $A \in \binom{[r]}{\leq k}$, we write $K_A(X)$
for the complete $|A |$-partite $|A|$-graph whose vertex classes are all the $X_i$ with $i\in A$.
The \emph{index} of an $r$-partite subset $S$ of $X$ is $i(S) = \{i \in [r]: S \cap X_i \neq \emptyset\}$.
Furthermore, given any set $B\subseteq i(S)$, we write $S_B = S \cap \bigcup_{i\in B} X_i $.
Similarly, given $A \in \binom{[r]}{\leq k}$ and an $r$-partite 
$k$-graph or $k$-complex $H$ on the vertex set $X$ we write $H_A$ for the collection of edges in $H$ of index $A$
and let $H_\emptyset=\{\emptyset\}$. 
In particular, if $H$ is a $k$-complex then $H_{\{i\}}$ is the set of all those vertices in~$X_i$
which lie in an edge of~$H$ (and thus form a (singleton) edge of~$H$).
In general, we will often view $H_A$ as an $r$-partite $|A|$-graph
with vertex set~$X$.  
Also, given a $k$-complex~$H$ we similarly write $H_{A^\leq} = \bigcup_{B \subseteq A} H_B$ and 
$H_{A^<} = \bigcup_{B \subset A} H_B$. We write $H^*_A$ for the $|A|$-graph whose edges are
those $r$-partite sets $S\subseteq X$ of index $A$ for which all proper subsets of $S$ belong to $H$. 
(In other words, a set $S$ with index $A$ satisfies 
$S \in H^*_A$ if and only if for all $j < |A|$ the edges of $H$ 
which have size $j$ and are subsets of $S$ form a complete $j$-graph on $|S|$ vertices.) 
Then the \emph{relative density of $H$ at index $A$} is 
$d_A(H) = |H_A|/|H_A^*|$. The \emph{absolute density of $H_A$} 
is $d(H_A) = |H_A|/|K_A(X)|$. (Note that  $|K_A(X)|=\prod_{i\in A} |X_i|$.)
If $H$ is a $k$-partite $k$-complex we may simply write $d(H)$ for $d(H_{[k]})$.
Similarly, the \emph{density} of a $k$-partite $k$-graph $H$ on $X=X_1 \cup \dots \cup X_k$
is $d(H)=|H|/|K_{[k]}(X)|$.

Finally, for any vertex $v$ of a hypergraph $H$, we define the \emph{vertex degree} $d(v)$ of $v$ to be the number 
of edges of $H$ which contain $v$. Note that this is not the same as the degree defined earlier, which was 
for sets of $k-1$ vertices. The \emph{maximum vertex degree} of $H$ is then the maximum of $d(v)$ taken 
over all vertices $v \in V(H)$. The \emph{vertex neighbourhood}  $VN(v)$ of $v$ is the set of all
vertices $u \in V(H)$ for which there is an edge of $H$ containing both $u$ and~$v$. For a 
$k$-partite $k$-complex $H$ on the vertex set $X_1 \cup \dots \cup X_k$ we also define the 
\emph{neighbourhood complex} 
$H(v)$ of a vertex $v \in X_i$ for some $i$ to be the $(k-1)$-partite $(k-1)$-complex with vertex set 
$\bigcup_{j\neq i} X_j$ and edge set $\{e \in H: e \cup \{x\} \in H\}$.

\subsection{Regular complexes}\label{regcomplexes}

In this subsection we shall define the concept of regular complexes (which was first introduced in the
$k$-uniform case by R\"odl and Skokan~\cite{RSk}) in the form used by R\"odl and Schacht~\cite{RS, RS2}.
This is a generalization of the 
standard concept of regularity in graphs, where we say that a bipartite graph $B$ on vertex classes $U$ and $V$ forms 
an $\eps$-regular pair if for any $U' \subseteq U$ and $V' \subseteq V$ with $|U'| > \eps |U|$ and $|V'| > \eps |V|$ 
we have $d(B[U' \cup V']) = d(B)\pm \eps$.

In the same way, we say that a $k$-complex $G$ is regular if the restriction of $G$ to any large subcomplex of 
lower rank has similar densities to $G$. More precisely, let $G$ be an $r$-partite $k$-complex on the vertex set 
$X = X_1 \cup \dots \cup X_r$. For any $A \in \binom{[r]}{\leq k}$, we say that $G_A$ is \emph{$\eps$-regular} if for 
any $H\subseteq G_{A^<}$ with $|H^*_A| \geq \eps |G^*_A|$ we have%
    \COMMENT{the definition below means we don't have to introduce restrictions to subcomplexes.
It seems to make sense if $|A|=1$}
$$
\frac{|G_A \cap H_A^*|}{|H_A^*|} = d_A(G) \pm \eps.
$$ 
We say $G$ is \emph{$\eps$-regular} if $G_A$ is $\eps$-regular for every $A \in \binom{[r]}{\leq k}$.
Note that if $G$ is a graph without isolated vertices, then the definition in the previous paragraph
is equivalent to the 2-complex $G^\le$ being $\eps$-regular.%
   \COMMENT{excluding isolated vertices in this remark is necessary, as isolated vertices of $G$
are not included as (singleton) edges in $G^\le$}
To illustrate the definition for $k=3$, suppose that $A=[3]$.
Then for instance the top level of $G_{[2]}$ is the bipartite subgraph of~$G$ induced by $X_1$ and $X_2$ and
$G_A^*$ is the set of (graph) triangles in $G$.
So roughly speaking, the regularity condition states that if we consider 
a subgraph of $G_{[2]} \cup G_{\{1,3\}} \cup G_{\{2,3\}}$
which spans a large number of triangles, 
then the proportion of these which also form an edge of $G_A$ is close to $d_A(G)$,
i.e.~close to the proportion of (graph) triangles in~$G$ between $X_1$, $X_2$ and $X_3$ which form an edge of~$G$.

Roughly speaking, the hypergraph regularity lemma states that 
an arbitrary $k$-graph can be split into pieces, each of which forms a regular $k$-complex. The version of 
the regularity lemma we shall use also involves the notion of a `partition complex', which is a certain partition 
of the edges of a complete $k$-complex.
As before, let $X = X_1 \cup \dots \cup X_r$ be an $r$-partite vertex set. A \emph{partition $k$-system~$P$ on $X$} 
consists of a partition $P_A$ of the edges of $K_A(X)$ for each $A \in \binom{[r]}{\leq k}$. We refer to 
the partition classes of $P_A$ as \emph{cells}. So every edge of $K_A(X)$ is contained in precisely
one cell of $P_A$. $P$ is a \emph{partition $k$-complex on $X$} if it also has the property that
whenever $S, S' \in K_A(X)$ lie in the same cell of $P_A$, we have that 
$S_B$ and $S_B'$ lie in the same cell of $P_B$ for any $B \subseteq A$. This property of $S, S'$ forms 
an equivalence relation on the edges of $K_A(X)$, which we refer to as \emph{strong equivalence}.
To illustrate this, again suppose that $k=3$ and $A=[3]$. Then if $P$ is a partition $k$-complex,
$P_{\{1\}}$, $P_{\{2\}}$ and $P_{\{3\}}$ together yield a vertex partition $Q_1$ refining
$X_1,X_2,X_3$. $Q_1$ naturally induces a partition $Q_2$
of the $3$ complete bipartite graphs induced by the pairs $X_i,X_j$. 
$P_{\{1,2\}}$, $P_{\{2,3\}}$ and $P_{\{1,3\}}$ also yield a partition $Q'_2$ of these
complete bipartite graphs. The requirement of strong equivalence now implies that $Q'_2$
is a refinement of $Q_2$. At the next level, $Q'_2$ naturally induces
a partition $Q_3$ of the set of triples induced by $X_1,X_2$ and $X_3$.
As before, strong equivalence implies that the partition $P_{\{1,2,3\}}$ of these triples
is a refinement of $Q_3$.

Let $P$ be a partition $k$-complex on $X=X_1 \cup \dots \cup X_r$. 
For $i\in [k]$, the cells of $P_{ \{i \} }$ are called \emph{clusters}
(so each cluster is a subset of some $X_i$).
We say that $P$ is \emph{vertex-equitable} if all clusters have the same size. 
$P$ is \emph{$a$-bounded} if $|P_A| \leq a$ for every $A$ (i.e.~if $K_A(X)$ is divided into at most $a$ cells by 
the partition $P_A$). Also, for any $r$-partite set%
    \COMMENT{use $Q$ instead of $S$ here since later on $S$ is used for something else in a similar
context, which we found confusing}
$Q\in \binom{X}{\le k}$, we write $C_Q$ for the set of
all edges lying in the same cell of $P$ as~$Q$, and write $C_{Q^\leq}$ for the $r$-partite $k$-complex
whose vertex set is $X$ and whose edge set is $\bigcup_{Q' \subseteq Q} C_{Q'}$. (Since~$P$ is a
partition $k$-complex, $C_{Q^\le }$ is indeed a complex.) The partition $k$-complex $P$ is 
\emph{$\eps$-regular} if $C_{Q^\leq}$ is $\eps$-regular for every $r$-partite $Q \in \binom{X}{\leq k}$. 

Given a partition $(k-1)$-complex~$P$ on $X$ and $A \in \binom{[r]}{k}$, we can define an 
equivalence relation on the edges of $K_A(X)$, namely that $S, S' \in K_A(X)$ are equivalent if and only if $S_B$ 
and $S_B'$ lie in the same cell of~$P$ for any strict subset $B \subset A$. We refer to this as
\emph{weak equivalence}. Note that if the partition complex $P$ is $a$-bounded, then $K_A(X)$ is
divided into at most $a^k$ classes by weak equivalence%
     \COMMENT{it really should be $a^k$, not $a$}. 
If we let $G$ be an $r$-partite $k$-graph on $X$, then we can use weak equivalence to refine the partition 
$\{G_A, K_A(X) \sm G_A\}$ of $K_A(X)$ (i.e.~two edges of $G_A$ are in the same cell if they are weakly 
equivalent and similarly for the edges not in $G_A$). 
Together with $P$, this yields a partition $k$-complex which we denote by $G[P]$. 
If $G[P]$ is $\eps$-regular then we say that $G$ is \emph{perfectly $\eps$-regular with respect to $P$}.
Note that if $G[P]$ is $\eps$-regular then $P$ must be $\eps$-regular too.

Finally, we say that $r$-partite $k$-graphs $G$ and $H$ on $X$ are \emph{$\nu$-close} if 
$|G_A \triangle H_A| < \nu |K_A(X)|$ 
for every $A \in \binom{[r]}{k}$, that is, if there are few edges contained in $G$ but not in $H$ and vice versa.

We can now present the version of the
regularity lemma we shall use to split our $k$-graph $H$ into regular $k$-complexes. It
actually states that there is some $k$-graph $G$ which is close to $H$ and which 
is regular with respect to some partition complex. This will 
be sufficient for our purposes, as we shall avoid the use of any edges in $G \sm H$, so every edge used will lie in 
both $G$ and $H$. There are various other forms of the regularity lemma for $k$-graphs
which give information on $H$ itself (the first of these were proved in~\cite{RSk,G1}) but 
these do not have the hierarchy of densities necessary for the application of the blow-up lemma (see~\cite{K2} for a 
fuller discussion of this point). The version below is due to R\"odl and Schacht~\cite{RS}
(actually it is a very slight restatement of their result).%
     \COMMENT{Weakened Theorem~\ref{eq-partition} to having just vertex-equitable since we didn't
use that \emph{all} the cells have the same size and now it's easier to believe that the Thm
follows from the (pf of the) R\"odl and Schacht result. Need that $a!r$ divides $n$ (instead of $a!$)
since we need that $a!$ divides each $|X_i|$. Also, removed that $P$ is $r$-partite, since by def all
our partition complexes on~$X$ are.}    

\begin{theo}[Theorem~14,~\cite{RS}]\label{eq-partition} 
Suppose integers $n,a,r,k$ and reals $\eps, \nu$ satisfy $1/n \ll\eps \ll 1/a \ll \nu, 1/r, 1/k$
and where $a!r$ divides~$n$. Suppose also that $H$ is an 
$r$-partite $k$-graph whose vertex classes $X_1, \dots, X_r$ form an equitable partition 
of its vertex set $X$, where $|X|=n$. Then there is an $a$-bounded $\eps$-regular
vertex-equitable partition $(k-1)$-complex $P$ on $X$ and an $r$-partite 
$k$-graph $G$ on $X$ that is $\nu$-close to $H$ and perfectly $\eps$-regular with respect to $P$.
\end{theo}

Here (and later on) we write $0<a_1 \ll a_2 \ll a_3 \ll a_4\leq  1$ to mean that we can choose the constants
$a_1,\dots,a_4$ from right to left. More
precisely, there are increasing functions $f_1,f_2,f_3$ such that, given
$a_4$, whenever we choose some $a_3 \leq f_3(a_4)$, $a_2\le f_2(a_3)$ and $a_1 \leq f_1(a_2)$, all
calculations needed in the proof of the subsequent statement are valid.
Hierarchies with more constants are defined similarly. 

One important property of regular complexes is that they remain regular when restricted to a large subset of their 
vertex set. For regular $k$-partite $k$-complexes this property is formalised by the following lemma, a special case of Lemma~6.18 in~\cite{K2}.

\begin{lemma}[Restriction of regular complexes]\label{restrict}
Suppose $\eps \ll \eps' \ll d  \ll c \ll 1/k$, and that $G$ is an $\eps$-regular $k$-partite 
$k$-complex on the vertex set $X = X_1 \cup \cdots \cup X_k$ such that $G_{\{i\}} = X_i$ for each $i$ and $d(G) > d$. 
Let $W$ be a subset of $X$ such that $|W\cap X_i| \geq c|X_i|$ for each $i$. 
Then the restriction $G[W]$ of $G$ to $W$ is $\eps'$-regular, with $d(G[W]) > d(G)/2$ and
$d_{[k]}(G[W]) > d_{[k]}(G)/2$.
\end{lemma}


\subsection{Robustly universal complexes.}

Apart from Theorem~\ref{eq-partition}, the other main tool
we shall use in the proof of Theorem~\ref{main} is the recent hypergraph blow-up 
lemma of Keevash. This result involves not only a $k$-complex $G$, 
but also a $k$-graph $M$ of `marked' edges on the same vertex set. If the pair $(G,M)$ is `super-regular', 
then this blow-up lemma can be applied to embed any spanning bounded-degree $k$-complex in $G \sm M$, that is, within 
$G$ but avoiding any marked edges. We will apply this with $M=G \setminus H$ where $G$ is the $k$-graph given
by Theorem~\ref{eq-partition}.
Super-regularity is a stronger notion than regularity. A result in~\cite{K2} states that every $\eps$-regular
$k$-complex can be made super-regular by deleting a few of its vertices.
Unfortunately, the notion of hypergraph super-regularity is very technical, but the following definition from~\cite{K2} avoids many of these technicalities. 
Let~$J'$ be a $k$-partite $k$-complex.
Roughly speaking, we say that $J'$ is robustly $D$-universal if the following holds:
even after the deletion of many vertices of $J'$, the resulting complex $J$ has the property that one can find in $J$ a copy of any $k$-partite $k$-complex $L$ which has vertex degree at most $D$ and whose vertex classes are the same as those of $J$. 
Condition (i) puts a natural restriction on the number of vertices we are allowed to delete from the 
neighbourhood complex of a vertex of $J$ and condition (iii) states that for a few vertices $u$ of $L$
we can even prescribe a `target set' in $V(J)$ into which $u$ will be embedded.

\medskip

\nib{Definition. (Robustly universal complexes)} Suppose that $J'$ is a $k$-partite $k$-complex on
$V' = V'_1 \cup\dots \cup V'_k$ with $J'_{\{i\}} = V'_i$ for each $i \in [k]$.  We say that $J'$ is \emph{$(c,c_0)$-robustly
$D$-universal} if whenever
\begin{itemize}
\item[(i)] $V_j \sub V'_j$ are sets with $|V_j| \ge c|V'_j|$ for all $j \in [k]$, such that writing $V = \bigcup_{j \in [k]} V_j$ and $J=J'[V]$ we have $|J(v)_=| \ge c|J'(v)_=|$ for any $j \in [k]$ and $v \in V_j$,
\item[(ii)] $L$ is a $k$-partite $k$-complex of maximum vertex degree at most $D$ on some vertex set $U = U_1 \cup \dots \cup U_k$ with $|U_j|=|V_j|$ for all $j \in [k]$,
\item[(iii)] $U_* \subseteq U$ satisfies $|U_* \cap U_j| \le c_0|U_j|$ for every $j \in [k]$, and sets $Z_u \subseteq V_{i(u)}$ satisfy $|Z_u| \ge c|V_{i(u)}|$ for each $u \in U_*$, where for each $u$ we let $i(u)$ be such that $u \in U_{i(u)}$,
\end{itemize}
\noindent then $J$ contains a copy of $L$, in which for each
$j \in [k]$ the vertices of $U_j$ correspond to the vertices of $V_j$, and
$u$ corresponds to a vertex of $Z_u$ for every $u \in U_*$.
\medskip

So our use of the blow-up lemma will be hidden through this definition. Of course, we shall also need to obtain 
robustly universal complexes. This is the purpose of the next theorem, which states that given a regular
$k$-partite $k$-complex $G$ with sufficient density, and a $k$-partite $k$-graph $M$ on the same vertex set
which is small relative to $G$, we can 
delete a small number of vertices from their common vertex set so that $G \sm M$ is robustly universal.
It is a special case of Theorem~6.32 in~\cite{K2}.

\begin{theo}\label{robust-universal}
Suppose that $1/n \ll \eps \ll c_0 \ll d^* \ll d_a \ll \theta \ll d, c, 1/k, 1/D, 1/C$,
$G$ is a $k$-partite $k$-complex on $V = V_1 \cup \dots \cup V_k$ 
with $n \leq |G_{\{j\}}| = |V_j| \leq Cn$ for every $j \in [k]$, 
$G$ is $\eps$-regular with $d_{[k]}(G) \ge d$ and $d(G_{[k]}) \ge d_a$, 
and $M \sub G_=$ with $|M| \le \theta |G_=|$. 
Then we can delete at most $2\theta^{1/3} |V_j|$ vertices
from each $V_j$ to obtain $V' = V_1' \cup  \dots \cup V_k'$, $G' = G[V']$ and $M' = M[V']$ 
such that
\begin{itemize}
\item[(i)] $d(G') > d^*$ and $|G'(v)_=| > d^*|G'_=|/|V'_i|$ for every $v \in V'_i$, and
\item[(ii)] $G' \sm M'$ is $(c,c_0)$-robustly $D$-universal.
\end{itemize}
\end{theo}

\section{Preliminary results} \label{sec:prelims}

In this section we will collect the preliminary results we need to prove Theorem~\ref{main}.
In order to apply Theorem~\ref{robust-universal}, we need to know under what conditions we can find particular
loose paths in complete $k$-partite $k$-graphs, which is the topic of the next subsection.

\subsection{Loose paths in complete graphs}

The problem of when we can find particular loose paths in a complete $k$-partite $k$-graph
can be reformulated in terms of the question of which strings satisfying certain adjacency conditions
can be produced from a fixed character set; the following lemma is the result we will need.

\begin{lemma} \label{strings}
Let $\ell$ and $a_1,\dots,a_k$ be integers such that
$0 \leq a_i < \ell/2$ for all $i$, and $\ell = \sum_{i=1}^k a_i$. Then for any
$s,t \in [k]$ there exists a string of length $\ell$ on alphabet $x_1,\dots,x_k$ such that
the following properties hold:
\begin{enumerate}
    \item[(1)] no two consecutive characters are equal,
    \item[(2)] the first character is not $x_s$ and the final character is not $x_t$,
    \item[(3)] the number of occurrences of character $x_i$ is $a_i$.
\end{enumerate} 
\end{lemma}
\nib{Proof.}
Note that the conditions on $\ell$ and the $a_i$ imply that $\ell\ge 3$.
We will construct the required string by starting with an `empty string' of~$\ell$
blank positions, and for each $i$ inserting precisely $a_i$ copies of character~$x_i$.
This ensures that condition~(3) will be satisfied. We shall fill the empty positions in the following order:
first the first position, then the third, and so on through the odd-numbered positions, until we reach
either position $\ell$ or position $\ell-1$ (dependent on whether $\ell$ is odd or even).
We then fill the second position, then the fourth, and so on until all positions are filled. Note
that if we proceed by inserting all copies of one character, then all the copies of another character,
and so forth, then condition~(1) must be satisfied. This is because to get two consecutive copies of
$x_i$, we must have inserted a copy of $x_i$ at some odd position $p$, then $p+2$, $p+4$, and
so on until reaching $\ell$ or $\ell-1$, and then filled even positions $2,4,6,\dots,p-1$. However,
this would imply that we had inserted at least $\ell/2$ copies of character $x_i$, contradicting
the fact that $a_i < \ell/2$.

We therefore only need to determine an order to insert the different characters so as to satisfy~(2).
We first consider the case $s \neq t$, say $s=1$ and $t=2$. In this case we insert $x_2$
first, $x_1$ last, and the remaining character blocks in any order in between. Clearly this prevents the
first character from being $x_1$ and the last from being $x_2$, and so~(2) is satisfied.
Now we may assume $s=t$, say $s=t=1$. Then if $\ell$ is odd, we insert the characters in the following order:
$x_2, x_3, \dots, x_k, x_1$. Then all the copies of $x_1$ must be in even positions
(since $a_1 < \ell/2$), and so (2) is satisfied. Alternatively, if $\ell$ is even, we
insert first $x_i$ for some $i \neq 1$ with $a_i > 0$, then $x_1$, and then the remaining blocks
of characters in any order. (Note that these include at least one character other than $x_1$ and $x_i$ since 
$\ell \ge 3$ and $a_j < \ell/2$ imply that at least three $j$ have $a_j \geq 1$.)
So neither the first nor last character can be $x_1$,
and so~(2) is again satisfied. \qed

\medskip

The next lemma is the result we were aiming for in this section, giving information about which
loose paths can be found in complete $k$-partite $k$-graphs. Note that the maximum vertex degree of a
loose path is two, and so this lemma will tell us when we can 
find a loose path in a robustly universal $k$-complex.
 
\begin{lemma}\label{loosepath}
Let $G$ be a complete $k$-partite $k$-graph on the vertex set $V_1\cup \dots\cup V_k$. Let $b_1,\dots,b_k$ be
integers with  $0\le b_i \leq |V_i|$ for each $i$. Suppose that
\begin{itemize}
    \item $n := \frac{1}{k-1}((\sum_{i=1}^{k} b_i) -1)$ is an integer, and
    \item $\frac{n}{2}+1 \leq b_i \leq n$ for all $i$.
\end{itemize}
Then for any $s,t \in [k]$, there exists a loose path in $G$ with an initial vertex in $V_s$,
a final vertex in $V_t$, and containing $b_i$ vertices from $V_i$ for each $i\in [k]$.
\end{lemma}
\nib{Proof.}
Note first that $n$ is the number of edges such a path must contain. Let $a_i = n - b_i$ for each $i$,
so that $0 \leq a_i < (n-1)/2$. By Lemma~\ref{strings} we can find a string
$S$ of length $n-1$ on the alphabet $V_1, V_2,\dots, V_k$ such that $V_i$ appears $a_i$ times,
no two consecutive characters are identical, the first character is not $V_s$ and the final character
is not $V_t$. Let $S_i$ be the $i$th character of $S$. To construct a loose path $P$ in $G$,
first choose any vertex from $V_s$ to be the initial vertex of $P$, and any vertex from $V_t$ to be the
final vertex of $P$. We also use $S$ to choose the link vertices of $P$: choose the $i$th link vertex
(i.e.~the vertex lying in the intersection of the $i$th and $(i+1)$th edges of $P$)
to be any member of $S_i$ not yet chosen. We have now assigned two vertices to each edge of~$P$.
Finally, we complete $P$ by assigning to each edge one as yet unchosen vertex from each of
the $k-2$ classes not yet represented in that edge. This is possible since precisely $a_i$ link
vertices are from the class~$V_i$ and so the total number of vertices used from
$V_i$ is $n-a_i = b_i$. Since $G$ is complete we know that each edge of~$P$ is an edge of $G$,
and so $P$ is a loose path satisfying all the conditions of the lemma.\qed

\subsection{Walks and connectedness in $k$-graphs} \label{walks}

A \emph{walk $W$} in a hypergraph $H$ consists of a sequence of edges $e_1, \dots, e_\ell$ of $H$ and
a sequence $x_0, \dots, x_{\ell}$ of (not necessarily distinct) vertices of $H$,
satisfying $x_{i-1} \neq x_i$ for all $i \in [\ell]$, and also $x_0 \in e_1$, $x_\ell \in e_\ell$
and $x_i \in e_i \cap e_{i+1}$ for all $i \in [\ell-1]$. The \emph{length of $W$} is the number of
its edges. We say that $x_0$ is the \emph{initial vertex of $W$}, $x_\ell$ is the
\emph{final vertex of $W$}, and that $x_1,\dots, x_{\ell-1}$ are the \emph{link vertices of $W$}.
By a \emph{walk from~$x$ to~$y$} we mean a walk with initial vertex~$x$ and final vertex~$y$.

Note that the vertices of a hypergraph~$H$ can be partitioned using the equivalence relation $\sim$,
where $x \sim y$ if and only if either $x=y$ or there exists a walk from $x$ to $y$.%
    \COMMENT{Reflexivity is given by the $x=y$ clause, symmetricity by the fact that we can reverse the walk
(i.e. take the $E_i$ and $x_i$ in the other order), and transitivity by noting we can concatenate
walks by adding the edges of the second walk on to the end of those of the first, and doing the same with
all but the initial $x_i$ of the second walk.}
We call the equivalence classes of this relation
\emph{components} of~$H$. We say that~$H$ is \emph{connected} if it has precisely
one component. Observe that all vertices of an edge of~$H$ must lie in the same 
component.%
    \COMMENT{Otherwise, this single edge is a walk from a vertex in one component to a vertex in
another, which is a contradiction}
Finally, note that if~$H$ is a connected hypergraph of order~$n$, then for any two vertices
$x,y$ of~$H$ we can find in a walk from~$x$ to~$y$ of length at most~$n$ in~$H$.%
     \COMMENT{Take $X_t = \{$ vertices which can be reached from a particular vertex
$v$ by a walk of length $ \leq t\}$, and $Y_t$ to be the complement of $X$. By connectivity there
must be an edge with a vertex in $X_t$ and a vertex in $Y_t$, and so $|X_{t+1}| > |X_{t}|$
unless $|X_{t}| = n$}

\subsection{Random splitting}

In this section we shall obtain, with high probability, a lower bound on the density of a subgraph of a
$k$-partite $k$-graph chosen uniformly at random. We will use Azuma's inequality on the
deviation of a martingale from its mean.

\begin{lemma}[Azuma~\cite{AZUMA}] \label{azuma}
Suppose $Z_0,\dots, Z_m$ is a martingale, i.e.~a sequence of random variables satisfying
$\mathbb{E}(Z_{i+1} \mid Z_0, \dots, Z_i) = Z_i$, and that $|Z_i-Z_{i-1}| \leq c_i$
for some constants $c_i$ and all $i\in [m]$. Then for any $t \geq 0$,
$$
\mathbb{P}(|Z_m - Z_0| \geq t) \leq 2 \exp \left( -\frac{t^2}{2 \sum_{i=1}^{m} c_i^2}\right).
$$
\end{lemma}

\begin{lemma} \label{randomsplit}
Suppose $1/n \ll c,\beta, 1/k, 1/b < 1$, and that $H$ is a $k$-partite $k$-graph on the vertex
set $X = X_1 \cup \dots \cup X_k$, where $n \leq |X_i| \leq bn$ for each $i\in [k]$. Suppose also that
$H$ has density $d(H) \ge c$ and that for each $i$ we have $\beta |X_i| \leq t_i \leq |X_i|$. 
If we choose a subset $W_i \subseteq X_i$ with $|W_i| = t_i$ uniformly at random
and independently for each $i$, and let $W = W_1 \cup \dots \cup W_k$, then the probability that 
$H[W]$ has density $d(H[W]) > c/2$ is at least $1-1/n^2$.
Moreover, the same holds if we choose $W_i$ by including each vertex of $X_i$ independently with 
probability~$t_i/|X_i|$.

\end{lemma}
\nib{Proof.} Let $m = |X|$. To prove the first assertion, 
we obtain our subsets $W_i \subseteq X_i$ through the following
two-stage random process, independently for each $i$. First we assign the vertices of each
$X_i$ into sets $X_i^1$ and $X_i^2$ independently at random, with each vertex being assigned
to $X_i^1$ with probability $t_i/|X_i|$, and assigned to $X_i^2$ otherwise. Then, in the (highly probable)
event that we have $|X_i^1| \neq t_i$ we shall select uniformly at random a set of vertices to transfer
between $X_i^1$ and $X_i^2$ to obtain from $X_i^1$ the set $W_i$ with $|W_i| = t_i$. For each $i$, no
subset $W_i \subseteq X_i$ of size $t_i$ is more likely to result from this process than any other, so
we have chosen each $W_i$ uniformly at random. It remains to show that $H[W]$ is likely to have high
density. We do this by noting that $H[X^1]$ is likely to have high density (where
$X^1 = X^1_1 \cup \dots \cup X^1_k$) and that with high probability we will only need to transfer a small
number of vertices to form $W=W_1\cup\dots\cup W_k$, which can have only a limited effect on the density.

More precisely, let $x_1, \dots, x_m$ be an ordering of the vertices of $X$, and for each $i\in [m]$
let the random variable $Y_i$ take the value 1 if $x_i \in X^1$, and 0 otherwise. Recall that we write $|H|$
to denote the number of edges of a $k$-graph $H$. For all $i=0,\dots,m$ we now
define random variables $Z_i$ by $Z_i = \mathbb{E}(|H[X^1]| \mid  Y_1, \dots, Y_i)$.
Then the sequence $Z_0, \dots, Z_m$ is a martingale, $Z_m = |H[X^1]|$, and as we formed each $X_i^1$ by assigning vertices of $X_i$ independently at random into $X_i^1$ and
$X_i^2$, we have $Z_0 = \mathbb{E}(|H[X^1]|) \ge c \prod_{i=1}^k t_i$. Also, for any vertex $x_i$, let $f(i)$ be such that $x_i \in X_{f(i)}$
(i.e.~$f(i)$ is the index of $x_i$).%
     \COMMENT{using $i(x_i)$ would have been awkard}
Then $|Z_i - Z_{i-1}| \leq \prod_{j \neq f(i)} |X_{j}| \leq (bn)^{k-1}$ for all $i\in [m]$. Thus we can apply
Lemma~\ref{azuma} to obtain
$$
\mathbb{P}\left( |Z_m - Z_0| \geq \frac{c\prod_{i=1}^k t_i}{4} \right) \leq
2 \exp \left( -\frac{c^2 \prod_{i=1}^k t_i^2}{32m b^{2k-2} n^{2k-2}}\right) \le \frac{1}{n^3}.
$$
Therefore the event that $d(H[X^1]) > 3c/4$ has probability at least $1-1/n^3$. Also, by a standard Chernoff bound,
for each $i\in [k]$ the event that $|X_i^1| = t_i \pm |X_i|^{2/3}$ has probability at least $1-1/n^3$. Thus with
probability at least $1-1/n^2$ all of these events will happen. Now, if $|X_i^1| > t_i$, we choose a set of
$|X_i^1| - t_i$ vertices of $X_i^1$ uniformly at random and move these vertices from $X_i^1$ to $X_i^2$.
Similarly, if $|X_i^1| < t_i$, then we choose a set of $t_i - |X_i^1|$ vertices of $X_i^2$ uniformly
at random and move these vertices to $X_i^1$. In either case, for any $i$ this action can decrease
$d(H[X^1])$ by at most $||X_i^1|-t_i|/|X_i^1| \ll c $. Thus if we let $W$ be the set obtained
from $X^1$ in this way, we have $d(H[W]) > c/2$, proving the first part of the lemma.

The proof of the `moreover part' is the same except that we can omit the `transfer' step at the end of the 
proof. \qed

\subsection{Decomposition of $G$ into copies of $\mc{A}_k$}
Let $\mc{A}_k$ denote the $k$-graph whose vertex set $V(\mc{A}_k)$ is the union of
$2k-2$ disjoint sets $U_0, U_1, U_2, \dots, U_{2k-3}$ of size $k-1$ and whose edges consist of all
$k$-tuples of the form $U_i \cup \{x\}$, with $i>0$ and $x \in U_0$ (see Figure~2).
\begin{figure}\label{fig:a_kdiag}
\centering\footnotesize
\psfrag{1}{$U_1$}\psfrag{2}{$U_2$}\psfrag{3}{$U_3$}\psfrag{4}{$U_0$}
\includegraphics[width=0.3\columnwidth]{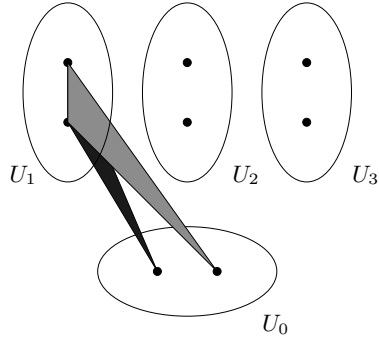}
\caption{The $3$-graph $\mc{A}_3$ (only edges involving $U_1$ are shown)}
\end{figure}
So $|V(\mc{A}_k)|=2(k-1)^2$. An \emph{$\mc{A}_k$-packing} in a $k$-graph~$G$ is a
collection of pairwise vertex-disjoint copies of $\mc{A}_k$ in~$G$. 
\begin{lemma} \label{ak-pack-lem}
Suppose $1/m \ll \theta \ll \psi\ll 1/k$, and that $G$ is a $k$-graph on $[m]$ such
that $|N_G(S)| > (\frac{1}{2(k-1)} + \theta)m$ for all but at most
$\theta m^{k-1}$ sets $S \in \binom{[m]}{k-1}$. Then~$G$ has an $\mc{A}_k$-packing which covers
more than $(1-\psi)m$ vertices of~$G$.
\end{lemma}
\nib{Proof.}
Let $A_1, \dots, A_t$ be an $\mc{A}_k$-packing of $G$ of maximum size, so $t \le m/(2(k-1)^2)$.
Let $X$ be the set of uncovered vertices, and suppose that $|X|>\psi m$.
Let $b = \theta |X|$. Our first aim is to choose disjoint sets $S_1,\dots,S_b$ in $\binom{X}{k-1}$
so that $|N_G(S_i)| > (1/(2(k-1)) + \theta)m$ and $|N_G(S_i) \cap X| < \theta m/2$ for all $i \in [b]$.
Note that  $\theta \ll \psi$ implies that
$\binom{|X|-2b(k-1)}{k-1} \gg\theta m^{k-1}$. So we can greedily choose disjoint
$S_1, \dots, S_{2b}\in \binom{X}{k-1}$
such that $|N_G(S_i)| > (1/(2(k-1)) + \theta )m$ for all $i\in [2b]$. 
Let $T = \{i \in [2b]: |N_G(S_i) \cap X| \geq \theta m/2\}$. We claim that $|T| \le b$.
Otherwise, consider the bipartite graph $B$ with vertex classes $T$ and $X$, where we join
$i \in T$ to $x \in X$ if $S_i \cup \{x\}$ is an edge of $G$. Note that $B$ cannot contain a
complete bipartite graph with $2k-3$ vertices in $T$ and $k-1$ vertices in $X$, as
this would correspond to a copy of $\mc{A}_k$ contained in~$X$, which is impossible as
$A_1, \dots, A_t$ is a maximum size $\mc{A}_k$-packing. However, by definition of ~$T$ we have
$d_B(i) \geq \theta m/2$ for every $i \in T$, and double-counting pairs $(i,P)$ with $i \in T$
and $P \in \binom{N_B(i)}{k-1}$ gives
$$|T| \binom{\theta m/2}{k-1} \leq \#\{(i,P)\} < (2k-3)\binom{|X|}{k-1},$$
a contradiction. This proves the claim, and by relabelling the $S_i$ we can assume that
$|N_G(S_i)| > (1/(2(k-1)) + \theta )m$ and $|N_G(S_i) \cap X| < \theta m/2$ for all $i \in [b]$.

Now we show how to enlarge the $\mc{A}_k$-packing~$A_1,\dots,A_t$. For $i \in [b]$ let
$$F_i = \{j\in [t]: |N_G(S_i) \cap V(A_j)| \ge k \}.$$
Since $|V(A_i)|= 2(k-1)^2$ for each $i\in [b]$ we have 
\begin{align*}
\left(\frac{1}{2(k-1)} + \frac{\theta}{2} \right)m
& <  |N_G(S_i) \sm X| = \sum_{j=1}^t |N_G(S_i) \cap V(A_j)|\\ & \le  |F_i| \cdot 2(k-1)^2 + (t-|F_i|) \cdot (k-1)
<  2(k-1)^2|F_i| +  \frac{(k-1)m}{2(k-1)^2},
\end{align*}
and so $|F_i| > \theta m/(4(k-1)^2)$. We now double-count pairs $(i,Q)$ with $i \in [b]$
and $Q \in \binom{F_i}{k-1}$. The number of such pairs is
$$\sum_{i=1}^b \binom{|F_i|}{k-1} > \theta \psi m \binom{\frac{\theta m}{4(k-1)^2}}{k-1} >
\sqrt{m}  \binom{t}{k-1}.$$
So we can find some $Q \in \binom{[t]}{k-1}$ and $R \subseteq [b]$ with $|R| > \sqrt{m}$
such that $Q \in \binom{F_r}{k-1}$ for every $r \in R$. For each $r \in R$ and each $q \in Q$
fix some $k$-set $K^{r,q} \subseteq N_G(S_r) \cap V(A_q)$ (which is possible by definition of $F_r$).
Then we can choose $R' \subseteq R$ with $|R'|=k(2k-3)$ so that $K^{r,q} = K^{r',q}$ for all
$r, r' \in R'$ and every $q \in Q$. For each $q \in Q$ we write $K^q$ for $K^{r,q}$ with $r \in R'$.

We will now use the~$K^q$ to find $k$ new copies of $\mc{A}_k$ that only intersect
$k-1$ of the copies in our packing. We arbitrarily divide $R'$ into $k$ sets
$R'_1, \dots, R'_k$ of size $2k-3$ and label $V(K^q) = \{v_{q,1}, \dots, v_{q,k}\}$
for all $q \in Q$. The new copies $A'_1,\dots,A'_k$ of $\mc{A}_k$ are obtained for
each $i \in [k]$ by identifying $U_1,\dots,U_{2k-3}$ with $\{S_r: r \in R'_i\}$ and
$U_0$ with $\{v_{q,i}\}_{q \in Q}$. Replacing the copies $\{A_q: q \in Q\}$ by $A'_1,\dots,A'_k$
we obtain a larger $\mc{A}_k$-packing. This contradiction completes the proof. \qed

\begin{coro} \label{ak-pack}
Lemma~\ref{ak-pack-lem} still holds if we insist that the sub-$k$-graph of $G$ induced by
the vertices covered by the $\mc{A}_k$-packing must be connected.
\end{coro}
\nib{Proof.}
Apply Lemma~\ref{ak-pack-lem} to obtain an $\mc{A}_k$-packing $A_1,\dots,A_\ell$ in $G$ with
$m_0:= |\bigcup_{i=1}^{\ell} V(A_i)| > (1- \psi/2)m$, and let $A$ be the
sub-$k$-graph of~$G$ induced by $\bigcup_{i=1}^{\ell} V(A_i)$.  By hypothesis at most $\theta m^{k-1}$
sets $S \in \binom{[m]}{k-1}$ have fewer than $m/(2(k-1))$ neighbours in $G$
and so at most $\theta m^{k-1}$ sets $T \in \binom{V(A)}{k-1}$ have no neighbours in~$V(A)$.
By the definition of a component, no edges of~$A$ contain vertices from different components of~$A$.
Therefore the largest component $C$ of $A$ must contain at least $(1-\psi)m$ vertices.
Indeed, if not then there are at least $\binom{m_0}{k-2}(\psi m/2)/(k-1)\gg \theta m^{k-1}$ sets
$T \in \binom{V(A)}{k-1}$ which meet at least two components of~$A$ and thus have no neighbours in~$A$,
a contradiction (we can obtain such a set $T$ by choosing $k-2$ vertices arbitrarily in $V(A)$ and then choosing
the final vertex in a different component of $A$ than the first vertex). 
Thus we may take the $\mc{A}_k$-packing consisting of all those copies~$A_i$
of $\mc{A}_k$ with $V(A_i)\subseteq V(C)$. \qed

\section{Proof of Theorem~\ref{main}}\label{proof}


In our proof we will use constants that satisfy the hierarchy
$$\frac{1}{n}  \ll \eps \ll d^*\ll d_a \ll \frac{1}{a} \ll \nu,\frac{1}{r} \ll \theta \ll d \ll c \ll \phi \ll \delta \ll \eta \ll \frac{1}{k}.$$ 
Furthermore, for any of these constants $\alpha$, we use $\alpha \ll \alpha' \ll \alpha '' \ll \dots$
and assume that the above hierarchy also extends to the additional constants, e.g.~$d''\ll c\ll c''\ll \phi$.

\subsection{Imposing structure on $H$}\label{structure}
\subsubsection{Step 1. Applying the regularity lemma}\label{regsection}
Let $H_1$ be the sub-$k$-graph obtained from~$H$ by removing up to $a!r$ vertices
so that $|V(H_1)|$ is divisible by $a!r$. Let $T = T_1 \cup \dots \cup T_r$ be an equitable $r$-partition of the
vertices of $H_1$, and let $H_2$ consist of all those edges of $H_1$ that are $r$-partite sets in $T$.
Then $H_2$ is an $r$-partite $k$-graph with order divisible by $a!r$, and so we may apply the regularity lemma
(Theorem~\ref{eq-partition}),
which yields an $a$-bounded $\eps$-regular vertex-equitable partition $(k-1)$-complex $P$ on $T$ and
an $r$-partite $k$-graph $G$ on $T$ that is $\nu$-close to $H_2$ and perfectly $\eps$-regular with respect to $P$. 

Let $M = G \sm H_2$. So any edge of $G\sm M$ is also an edge of $H$. Let
$V_1,\dots, V_m$ be the clusters of~$P$. So $T=V_1\cup \dots\cup V_m$ and $G$ is $m$-partite
with vertex classes $V_1\cup \dots\cup V_m$.%
     \COMMENT{the latter holds since $V_1,\dots,V_m$ refines $T_1,\dots,T_r$}
Note that $m \leq ar$ since $P$ is $a$-bounded. Moreover, since $P$ is vertex-equitable,
each $V_i$ has the same size. So let $n_1 = |V_i| =|T|/m$. 

As is usual in regularity arguments, we shall consider a reduced $k$-graph, whose vertices correspond to
the clusters $V_i$, and whose edges indicate that within the cells of~$P$ corresponding to the edge we
can find a subcomplex to which we can apply Theorem~\ref{robust-universal}. For this we would
like $G$ to have high density in these cells, and $M$ to have low density. Thus we define the \emph{reduced
$k$-graph $R$} on $[m]$ as follows: a $k$-tuple $S$ of vertices of $R$ corresponds to the $k$-partite
union $S' = \bigcup_{i \in S} V_i$ of clusters. The edges
of $R$ are precisely those $S\in \binom{[m]}{k}$ for which $G[S']$ has density at least $c''$
(i.e. $|G[S']| > c'' |K_S(S')|$) and for which $M[S']$ has density at most $\nu^{1/2}$ (i.e.
$|M[S']| < \nu^{1/2}|K_S(S')|$). 

Now, the edges in the reduced graph are useful in the following way.
Given an edge $S\in R$, let $S' = \bigcup_{i \in S} V_i$ again.
Using weak equivalence (defined in Section~\ref{regcomplexes}), the cells of~$P$ induce
a partition $C^{S,1},\dots,C^{S,m_S}$ of the edges of $K_S(S')$. Recall that $m_S\le a^k$.
Therefore at most $c''|K_S(S')|/3$ edges of $K_S(S')$ can lie in sets $C^{S,i}$
with $|C^{S,i}| \leq c''|K_S(S')|/(3a^k)$. Furthermore, $|M[S']| < \nu^{1/2} |K_S(S')|$
(as $S\in R$) and so at most $\nu^{1/4}|K_S(S')|$ edges of $K_S(S')$ can lie in sets $C^{S,i}$ with
$|M \cap C^{S,i}| \geq \nu^{1/4}|C^{S,i}|$. Together with the fact that $|G[S']| > c'' |K_S(S')|$ this
now implies that more than $c''|K_S(S')|/2$ edges of $G[S']$ lie in sets $C^{S,i}$ with
$|C^{S,i}| > c''|K_S(S')|/(3a^k)$ and $|M \cap C^{S,i}| < \nu^{1/4}|C^{S,i}|$.
Thus there must exist such a set $C^{S,i}$ that also satisfies $|G \cap C^{S,i}| >c''|C^{S,i}|/2$.
Fix such a choice of $C^{S,i}$ and denote it by $C^S$.
Let $G^S$ be the $k$-partite $k$-complex on the vertex set $S'$ consisting of $G \cap C^S$ and the
cells of $P$ that `underlie' $C^S$, i.e. for any edge $Q\in G\cap C^S$ we have
\begin{equation}\label{eqGS}
G^S=(G\cap C^S)\cup \bigcup_{Q'\subset Q} C_{Q'}.
\end{equation}
(Recall that $C_{Q'}$ was defined in Section~\ref{regcomplexes}.)
We also define the $k$-partite $k$-graph $M^S = G^S \cap M$ on the vertex set~$S'$.
Then the following properties hold:
\begin{itemize}
\item[(A1)] $G^S$ is $\eps$-regular. 
\item[(A2)] $G^S$ has $k$-th level relative density $d_{[k]}(G^S) \geq d'$. 
\item[(A3)] $G^S$ has absolute density $d(G^S) \geq d_a'$. 
\item[(A4)] $M^S$ satisfies $|M^S| < 2\nu^{1/4}|(G^S)_=|/c''$.
\item[(A5)] $(G^S)_{\{i\}} = V_i$ for any $i \in S$.
\end{itemize}
Indeed, (A1) follows from~(\ref{eqGS}) since $G$ is perfectly $\eps$-regular with respect to $P$.
To see~(A2), note that $(G^S_{[k]})^*=C^S$ and so 
$d_{[k]}(G^S) = |G^S_{[k]}|/|(G^S_{[k]})^*|=|G^S \cap C^S|/|C^S| > c''/2$
by our choice of $C^S$. Similarly, (A3) follows from our choice of~$C^S$ since
$$d(G^S) = \frac{|G^S_{[k]}|}{|K_S(S')|} = \frac{|G^S \cap C^S|}{|C^S|} \cdot \frac{|C^S|}{|K_S(S')|} > \frac{(c'')^2}{6a^k}>d'_a.$$ 
(A4) holds since $|(G^S)_=| = |G \cap C^S| >c''|C^S|/2$ and $|M^S| \leq |M \cap C^S| < \nu^{1/4}|C^S|$.
Finally, (A5) follows from~(\ref{eqGS}) and the fact that $C_{\{v\}}=V_i$ for all $v\in V_i$.

\subsubsection{Step 2. Choosing an $\mc{A}_k$-packing of~$R$}
The next step in our proof is to use Corollary~\ref{ak-pack} to find an $\mc{A}_k$-packing in the
reduced $k$-graph~$R$. For this we shall need an approximate minimum degree condition for~$R$. Let
 $$J = \left\{I \in \binom{[m]}{k-1}: |N_R(I)| \le
 \left(\frac{1}{2(k-1)} + \phi \right)m \right\}.$$
We shall show that $J$ is small, that is, that almost all $(k-1)$-tuples of vertices of~$R$ have degree at
least $(1/(2(k-1))+\phi)m$ in~$R$. Consider how many edges of $H$ do not belong to $G[S']$
for some edge $S \in R$. (Recall that $S'=\bigcup_{i\in S} V_i$.) There are three possible reasons why an
edge $e\in H$ does not belong to such a restriction:
\begin{itemize}
    \item[(i)] $e$ is not an edge of $G$. This could be because $e$ lies in $H$ but not $H_1$,
in $H_1$ but not $H_2$, or in $H_2$ but not $G$. There are at most $a!rn^{k-1}$ edges of the
first type, at most $n^k/r$ of the second type, and at most $\nu n^k$ of the third type.
    \item[(ii)] $e \in G$ contains vertices from $V_{i_1}, \dots, V_{i_k}$ such that the restriction of $M$ to
$S'=\bigcup_{i\in S} V_i$ satisfies $|M[S']| \geq \nu^{1/2}|K_S[S']|$, where $S=\{i_1,\dots,i_k\}$.
(Note that since $G$ and thus $M$ is $m$-partite, $i_1, \dots, i_k$ are all distinct.)
Since $G$ and~$H_2$ are $\nu$-close and thus $|M|\le \nu n^k$ there are at most $\nu^{1/2} n^k$ edges
of this type. 
    \item[(iii)] $e \in G$ contains vertices from $V_{i_1}, \dots, V_{i_k}$ such that the restriction of
$G$ to $\bigcup_{i\in S} V_i$ has density less than $c''$. There are at most $c''n^k$ edges of this type. 
\end{itemize}
Therefore there are fewer than $2c''n^k$ edges of $H$ that do not belong to the restriction of
$G$ to $S'$ for some $S \in R$, and so we have
\begin{align*}
|J| n_1^{k-1} \cdot \left( \frac{1}{2(k-1)} + \eta \right)n
& < \sum_{I \in J} \sum_{x_i \in V_i, i \in I} |N_H(\{x_i : i \in I\})| \\
 < 2c''k n^k + \sum_{I \in J} |N_R(I)|n_1^k
 &\leq  2c'' kn^k + |J|\left(\frac{1}{2(k-1)} + \phi \right)m n_1^k.
\end{align*}
Since $n-a!r\le mn_1\le  n$ we deduce that $|J|n_1^{k-1}(\eta-\phi)n < 2c''kn^k<3c''k(mn_1)^{k-1}n$, and so
$|J| < \phi m^{k-1}$
(since $c'' \ll \phi \ll \eta$). This allows us to apply Corollary~\ref{ak-pack} (with $G=R$) to
obtain an $\mc{A}_k$-packing $A_1, \dots, A_t$ in $R$ with $|\bigcup_{i=1}^t A_i| > (1-\delta)m$,
such that the sub-$k$-graph of $R$ induced by $\bigcup_{i=1}^t V(A_i)$ is connected. For each $i\in [t]$,
let the vertex set of $A_i$ be $U^i_0 \cup U^i_1 \cup \dots \cup U^i_{2k-3}$, with each $U^i_j$
of size $k-1$, so that the edge set is $\{ U^i_j \cup \{x\}: j\in [2k-3], x \in U^i_0\}$.

\subsubsection{Step 3. Forming the exceptional path.}
Given a sub-$k$-graph $R'$ of~$R$ and a cluster~$V_i$, we say that \emph{$V_i$ belongs to~$R'$}
if $i\in V(R')$. Let $V'_0$ contain the at most $a!r$ vertices of $H$ we removed at the start of the proof,
and also the vertices in all those clusters not belonging to some copy of $\mc{A}_k$
in our packing (there are at most $\delta n$ of the latter).
We will incorporate these vertices into a path~$L_e$ which will later form part of our
loose Hamilton cycle. 
We also include in $V'_0$ an arbitrary choice of $\delta n_1$ vertices from each $V_y$ for which
$y \in U_j^i$ for some $j\in [2k-3]$ and some $i\in [t]$ (we do not modify any of the $V_y$
for which $y \in U_0^i$). 
We add up to $k-3$ more vertices from $U^1_1$ (say)%
    \COMMENT{all we need is that we don't touch the $U^i_0$s}
to $V'_0$ so that%
    \COMMENT{previous method mixed up mod $k-1$ and mod $k-2$, so didn't work}
$|V'_0| \equiv 0 \mod k-2$. We delete all these vertices
from the clusters they belonged to and still write~$V_y$ for
the subcluster of a cluster~$V_y$ obtained in this way.
This gives $|V'_0| \leq 5\delta n/2$. 

Now, we shall construct a path $L_e$ in $H$, which will contain all the vertices in~$V_0'$ and
avoid all the clusters~$V_y$ with $y \in U_0^i$.  Let
$V_{>0} = \bigcup\{V_y: y \in U_j^i, j\in [2k-3], i\in [t]\}$.
So we shall use only vertices from $V'_0$ and $V_{>0}$ in forming $L_e$. Recall that
if $|V(H)|$ is not a multiple of $k-1$, then a loose Hamilton cycle contains a single pair of edges
which intersect in more than one vertex: we shall make allowance for this here. Choose $A,B \subseteq V_{>0}$
satisfying $|A|=|B| = k-1$, $|A \cap B| \equiv 1-|V(H)|\mod k-1$ and $1 \leq |A \cap B| \leq k-1$.
Now choose distinct $x_0, x_1 \in V_{>0}\setminus (A\cup B)$ such that $\{x_0\} \cup A \in H$ and
$\{x_1\} \cup B \in H$ (we shall see in a moment that such $x_0,x_1$ exist). These edges will be the first 2
edges of~$L_e$. To complete $L_e$, let $Z_1,\dots,Z_s$ be any partition of the vertices of
$V'_0$ into sets of size $k-2$. We proceed greedily in forming $L_e$: for each $i=1,\dots,s$ choose any $x_{i+1}\in V_{>0}\setminus (A\cup B)$ such that $Z_i \cup \{x_i,x_{i+1}\} \in H$
(where the $x_i$ are all chosen to be distinct). 

Let us now check that there will always be such a vertex available.
Indeed, every set in $\binom{V(H)}{k-1}$
has at least $(1/(2(k-1))+\eta)n$ neighbours and we can choose any such neighbour which lies in~$V_{>0}$
and has not already been used. But $|V(H) \sm V_{>0}| \leq n/(2(k-1)) +|V_0'|$ and at most
$|V_0'|+2k\le 3\delta n$ vertices have been used before. Thus (since $\delta \ll \eta$) for each choice of an~$x_i$
we have at least $\eta n/2$ vertices of $V_{>0}$ to choose from. Moreover, these vertices must be
contained in at least $\eta n/(2n_1)$ different $V_y$ such that $y \in U_j^{i'}$ ($j>0$).
Thus we can avoid choosing a vertex from any single~$V_y$ more than
$6 \delta n_1/\eta \le \delta' n_1/2$ times. The path $L_e$ thus formed has edges
$\{x_0\} \cup A$, $B \cup \{x_1\}$ and $\{x_i, x_{i+1}\} \cup Z_i$ for all $i \in [s]$.
So all the vertices of $V'_0$ are included in $L_e$. 
For each cluster $V_y$, we still denote the subset of~$V_y$ lying in $V(H-L_e)$ by $V_y$.
Then each $V_y$
with $y \in U_0^i$ for some $i$ still satisfies $|V_y| = n_1$, and each $V_y$ with $y \in U_j^i$ for
some $j>0$ satisfies
\begin{equation} \label{Vy}
(1-\delta')n_1 \leq \left(1-\delta-\frac{\delta'}{2}\right)n_1-(k-3)\le |V_y|\leq (1-\delta)n_1.
\end{equation}
In addition%
    \COMMENT{$|A\cup B|\equiv (2k-2)-1+|V(H)|\mod k-1$} 
\begin{equation} \label{exceppathsize}
|V(H) \sm V(L_e)| \equiv |V(H)|-|A \cup B \cup \{x_0,x_1 \}| \equiv -1 \mod k-1.
\end{equation} 
Note that $L_e$ need not be a loose path, but that even if it is not it may still form part of a
loose Hamilton cycle. Also observe that $|V(L_e)| \leq 6 \delta n$.

\subsubsection{Step 4. Splitting our copies of $\mc{A}_k$.}
The next step of the proof will be to split the copies~$A_1,\dots,A_t$ of $\mc{A}_k$
(more precisely the clusters belonging to the~$A_i$) into sub-$k$-complexes of $G$
that we shall later use to embed spanning loose paths. 
Consider any $A_i$. For convenient notation we identify each $U^i_j$ in $A_i$
with $[k-1]$ (but recall that they are disjoint sets). For each $y \in U^i_0 = [k-1]$ we
have $|V_y| = n_1$, and so we can partition $V_y$ uniformly at random
into $2k-3$ pairwise disjoint subsets $S^i_{y,1}, \dots, S^i_{y,2k-3}$, each of size%
    \COMMENT{by def of $n_1$, these are all integers}
$\frac{n_1}{2k-3}$. Similarly, given $z\in U^i_j=[k-1]$ with $j\in [2k-3]$,
(\ref{Vy}) and the  fact that $\delta' \ll \eta$ imply that we can partition
$V_z$ uniformly at random into $k-1$ pairwise disjoint subsets $T^i_{j,z}$ and
$\{U^i_{j,z,w}\}_{w \in [k-1] \sm \{z\}}$ so that
$\frac{n_1}{2k-3} \leq |T^i_{j,z}| \leq \frac{(1-\eta)2n_1}{2k-3}$ and
$|U^i_{j,z,w}| = \frac{(1-\eta)2n_1}{2k-3}$ for all $w \in [k-1] \sm \{z\}$.
Figure~3 shows how we do this in the case $k=3$.
\begin{figure}\label{fig:split}
\centering\footnotesize
\psfrag{1}{$U^i_1$} \psfrag{2}{$U^i_2$} \psfrag{3}{$U^i_3$} \psfrag{4}{$U^i_0$}
\psfrag{5}{$\! \! U^i_{1,2,1}$} \psfrag{6}{$T^i_{1,2}$} \psfrag{7}{$U^i_{1,1,2}$}
\psfrag{8}{$T^i_{1,1}$} \psfrag{9}{$U^i_{2,1,2}$} \psfrag{10}{$T^i_{2,1}$}
\psfrag{11}{$S^i_{1,1}$} \psfrag{12}{$S^i_{1,2}$} \psfrag{13}{$S^i_{1,3}$}
\psfrag{14}{$S^i_{2,1}$} \psfrag{15}{$S^i_{2,2}$} \psfrag{16}{$S^i_{2,3}$}
\psfrag{17}{$T^i_{2,2}$} \psfrag{18}{$U^i_{2,2,1}$} \psfrag{19}{$T^i_{3,1}$}
\psfrag{20}{$U^i_{3,1,2}$} \psfrag{21}{$T^i_{3,2}$} \psfrag{22}{$U^i_{3,2,1}$}
\includegraphics[width=0.4\columnwidth]{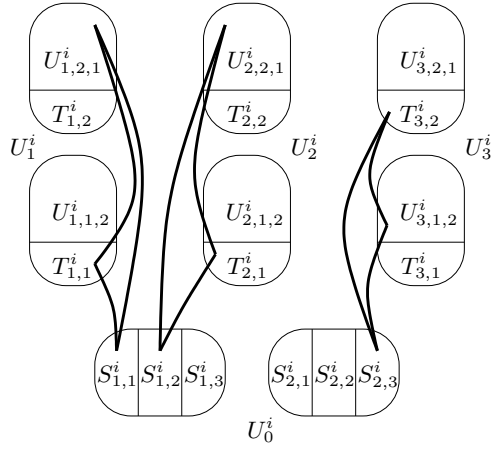}
\caption{Splitting up $A_i$ in the case $k=3$.}
\end{figure}

We arrange these pieces
into $(k-1)(2k-3)$ collections of $k$ sets as follows: for each $y \in U^i_0$ and each
$j \in [2k-3]$ we have a collection consisting of $S^i_{y,j}$, $T^i_{j,y}$ and $\{U^i_{j,z,y}\}_{z \ne y}$.
($3$~of these collections are illustrated in Figure~3.)
For convenient notation we relabel these collections as $\{ X_{i,1}, \dots, X_{i,k} \}$
with $1 \le i \le t' = (k-1)(2k-3)t$, where for all $i \in [t']$ we have
\begin{equation} \label{splitsizes1}
|X_{i,1}| = \frac{n_1}{2k-3},\  \frac{n_1}{2k-3} \leq |X_{i,2}| \leq 
\frac{(1-\eta)2n_1}{2k-3} \textrm{ and } |X_{i,j}| = \frac{(1-\eta)2n_1}{2k-3} 
\ \mbox{ for } 3 \le j \le k,
\end{equation}
and
\begin{equation} \label{splitsizes2}
(1-\delta')n_1 \le \sum_{j=2}^k  |X_{i,j}| \le (1-\delta)n_1 
\end{equation}
((\ref{splitsizes2}) follows from~(\ref{Vy}) using the fact that all the $U^{i'}_{j',z,w}$ have equal size.)
Let $X_i = \bigcup_{j \in [k]} X_{i,j}$,
so each $X_i$ is a $k$-partite set, on which we shall now find a sub-$k$-complex~$G_i$ of
$G$ that is suitable for applying Theorem~\ref{robust-universal}. 

Consider any copy $A_{i'}$ in our $\mc{A}_k$-packing. Note that for each of the $(k-1)(2k-3)$
collections $\{X_{i,1}, \dots, X_{i,k} \}$ obtained by splitting up the clusters belonging
to~$A_{i'}$ there is an edge $S(i)\in A_{i'}$ such that each~$X_{i,j}$ lies in a cluster belonging to~$S(i)$
(and these clusters are distinct for each of $X_{i,1},\dots,X_{i,k}$).
Recall that $S'(i)$ denotes the union $\bigcup_{\ell\in S(i)} V_\ell$ of all the clusters belonging
to~$S(i)$. Let $G_i$ denote the restriction of the $k$-partite $k$-complex $G^{S(i)}$ (which
was defined in Section~\ref{regsection}) to $X_i$, i.e.~$G_i=G^{S(i)}[X_i]$. Let $M_i = M \cap G_i=M^{S(i)}[X_i]$.
We claim that we may choose the above collections $\{ X_{i,1}, \dots, X_{i,k} \}$ such that
\begin{equation}\label{desityHG}
d(H[X_i])\ge \frac{c''}{4} \  \ \text{for all } i\in [t'].
\end{equation}
Indeed, since $S(i)\in R$, $G[S'(i)]$ has absolute density at least $c''$ and $M[S'(i)]$ has density
at most $\nu^{1/2}$. Since $G\sm M\subseteq H$ and $\nu \ll c''$ this shows that $H[S'(i)]$ has density at
least $c''/2$. Lemma~\ref{randomsplit}
now implies that each $H[X_i]$ has density at least $c''/4$ with probability $1-1/n_1^2$, and so with 
non-zero probability this is true for all $i\in [t']$. 

Lemma~\ref{restrict} and properties (A1)--(A3) and~(A5)
imply that~$G_i$ is an $\eps'$-regular $k$-partite $k$-complex on the vertex set $X_i$,
with absolute density $d(G_i) \geq d(G^{S(i)})/2\ge d_a$, relative density $d_{[k]}(G_i) \geq d$, and
$(G_i)_{\{j\}} = X_{i,j}$ for each $j$.
Moreover, using $\nu \ll \theta \ll c$, property~(A4) and the fact that $d(G_i) \geq d(G^{S(i)})/2$
we see that
$$|M_i| \le  |M^{S(i)}|<\frac{2\nu^{1/4}|(G^{S(i)})_=|}{c''}\le \theta |(G_i)_=|.$$
So by Theorem~\ref{robust-universal} we can delete at most $\theta' |X_{i,j}|$ vertices from
each $X_{i,j}$  so that if we let $X'_{i,j} \subseteq X_{i,j}$ consist of the undeleted vertices, and let $X'_i := \bigcup_{j=1}^k X'_{i,j}$, $G_i' := G_i[X'_i]$ and $M_i' := M_i[X'_i]$, then $G_i' \sm M_i'$ is $(c,\eps'')$-robustly $2^k$-universal, $d(G_i') > d^*$ and $|G_i'(v)_=| > d^*|(G_i')_=|/|X'_{i,j}|$ for every $v \in X'_{i,j}$. In particular, the latter two conditions together imply that $d(G_i'(v)_=) > (d^*)^2$ for every $v \in X_i'$. Let $X''$ denote the set of vertices deleted from any $X_{i,j}$, so $|X''|\le \theta' n$. By deleting up to $k-3$ more vertices if necessary, we may assume that $|X''|$ is divisible by~$k-2$. The latter will help us to extend~$L_e$ into a
path which contains all the vertices in~$X''$.

\subsubsection{Step 5. Extending the exceptional path~$L_e$.}
When extending~$L_e$ in order to incorporate~$X''$, we shall have to remove some more vertices from
some of the $X'_{i,j}$, and we wish to do this so that the remainder satisfies (i) in the 
definition of robust universality.
For this reason, we partition each $X'_{i,j}$ into
two parts $AX'_{i,j}$ and $BX'_{i,j}$ as follows (where we write $BX'_i$ for
$\bigcup_{j\in [k]} BX'_{i,j}$):
\begin{itemize}
    \item[(B1)] For all $i,j$ and every $v \in X'_{i,j}$ we have $|(G_i'(v)[ BX'_i])_= | \ge 2c|G_i'(v)_=|$. 
    \item[(B2)] Every set of $k-1$ vertices of $H$ has at least $n/(4k)$ neighbours in $\bigcup_{i,j} AX'_{i,j}$. 
\end{itemize}
(Recall that for a $(k-1)$-complex $F$, $F_=$ denotes the `$(k-1)$th level' of $F$.)
To see that such a partition exists, consider a partition obtained by assigning each vertex to a part
with probability $1/2$ independently of all other vertices.
(B2) is then satisfied with high probability by a standard Chernoff bound.
Now consider (B1). The `moreover' part of 
Lemma~\ref{randomsplit} implies that with high probability we have for all $i, j$ and for all
$v \in X'_{i,j}$ that $d((G_i'(v)[ BX'_i])_=) \ge d(G_i'(v)_=)/2$.
Also, a standard Chernoff bound implies that with high probability $|BX'_{i,j'}|\ge |X'_{i,j'}|/3$
for all $j'\in [k]$. Thus
$$ |(G_i'(v)[ BX'_i])_= |=d((G_i'(v)[ BX'_i])_=)\prod_{j'\neq j}|BX'_{i,j'}|\ge
\frac{d(G_i'(v)_=)}{2}\prod_{j'\neq j}\frac{|X'_{i,j'}|}{3}\ge 2c|G_i'(v)_=|.
$$
Now, we shall extend our path~$L_e$ to include the vertices in~$X''$, using only vertices
from $\bigcup_{i,j} AX'_{i,j}$. 
We proceed similarly to when constructing~$L_e$. So we split~$X''$ into sets $Z_1,...,Z_{s'}$
of size $k-2$ (so $s' \leq \theta' n$). Letting $x_0$ be a final vertex of $L_e$,
for $i\in [s']$, we successively choose $x_i$ to be a neighbour of the $(k-1)$-tuple $Z_i \cup \{x_{i-1}\}$ contained
in some $AX'_{i',j'}$ and not already included in $L_e$, and extend $L_e$ by the edge
$Z_i \cup \{x_{i-1}, x_i\}$, continuing to denote the extended path by~$L_e$.
Recall that $L_e$ originally contained at most $6 \delta n$ vertices. Since $|X''|\le \theta' n$,
after each extension of $L_e$ we shall have $|V(L_e)| < \eta n$. So~(B2) implies that for each choice of $x_i$
we have at least $n/(5k)$ suitable vertices and hence at least $t'/(5k)$ of the sets $AX'_{i'}$ contain
such a suitable vertex. This shows that we can choose the~$x_i$ in such a way that
at most $\theta'' n_1$ vertices are chosen from any single $AX'_{i'}$.

For each~$i\in [t']$ let $X^i = X^i_1 \cup \dots \cup X^i_k$ be the vertices
remaining after the removal from $X_i'$ of the at most $\theta'' n_1$ vertices used in extending $L_e$, let $G^i = G'_i[X^i]$, and let $M^i = M'_i[X^i]$.
By~(\ref{desityHG}) there are at least $cn$ vertices $v \in V(H)$ such that $v$ lies in some $X^i$ for which at least $|H[X^i]|/(2|X^i|)$ edges of $H[X^i]$ contain $v$. So we may add two further edges of $H$ to $L_e$ (one at each end) so that the new path $L_e$ has an initial vertex $x_e$ and a final vertex $y_e$ which each lie in at least $|H[X^i]|/(2|X^i|)$ edges of their respective $H[X^i]$. (We also delete the vertices of these additional two edges from their $X^i$, $G^i$ and $M^i$). Note that $x_e$ may be contained in some $BX'_{i,j}$ (and the same is true of $y_e$), but by (B2) we may choose these two additional edges so that all other vertices used lie in some $AX_{i,j}'$.



%
   \COMMENT{the average vertex degree is larger than $|H[X^i]|/(2|X^i|)$, so choose such a vertex as $x_e$ 
and use the mindeg property to connect $L_e$ to $x_e$. The factor 2 is only there to ensure that $x_e \neq y_e$,
though with more thought one can actually guarantee even without the factor 2...}

We claim that the above steps give us the following useful structure: 
a path $L_e$ which is ready to form part of a 
loose Hamilton cycle, and disjoint $k$-partite vertex sets $X^i = X^i_1 \cup \dots \cup X^i_k$ supporting $k$-complexes $G^i$ and $k$-graphs $M^i$ for each $i \in [t']$ which satisfy the following properties.
	\begin{itemize}
	\item[(C1)] Every vertex of $H$ lies in either the path $L_e$ or precisely one of the $k$-partite sets $X^i$.
	\item[(C2)] For each $i$, $G^i$ is a $k$-partite sub-$k$-complex of $G$ on the vertex set $X^i$. $M^i$ is the $
k$-partite $k$-graph $M \cap G^i$, and $G^i \sm M^i \subseteq H$.  Clearly these statements remain true after the
 deletion of up to $\eps n_1$ vertices of $X^i$.
	\item[(C3)] Even after the deletion of up to $\eps n_1$ vertices of $X^i$, the following statement holds. 
Let $L$ be a $k$-partite $k$-complex on the vertex set $U = U_1 \cup \dots \cup U_k$, where $|U_j| = |X^i_j|$ for 
each $j$, and let $L$ have maximum vertex degree at most $2^k$. Let $\ell \le 2(t')^2$ and
suppose we have $u_1, \dots, u_\ell \in U$ and sets 
$Z_s \subseteq X^i_{j(u_s)}$ with $|Z_s| \geq c|X^i_{j(u_s)}|$ for each $s \in [\ell]$ (where $j(u_s)$ is such that $u_s \in U_{j(u_s)}$). 
 Then $G^i \sm M^i$ contains a copy of $L$, in which for each $j$
the vertices of $U_j$ correspond to the vertices of $X^i_j$, and each $u_s$ corresponds to a 
vertex in $Z_s$.
   \item[(C4)] For each $i$, $H^i = H[X^i]$ has density at least $c'$, even after the deletion of up to 
$\eps n_1$ vertices of $X^i$.
   \item[(C5)] If we delete up to $\eps n_1$ vertices from any $X^i$, and let $t_j = |X^i_j|$ for each $j\in [k]$ after 
these deletions, and let $n'_i = \frac{(\sum t_j)-1}{k-1}$, then $n'_i/2+1 \leq t_j \leq n'_i$ for all $j$.
   \item[(C6)] The initial vertex $x_e$ of $L_e$ lies in at least $|H[X^i]|/(2|X^i|)$ edges of~$H[X^i]$,
where $i$ is such that $x_e\in X^i$. The analogue holds for the final vertex~$y_e$ of~$L_e$. 
	\end{itemize}
(When we talk of removing a vertex of $X^i$ we implicitly mean that $G^i$, $M^i$ and $H^i$ are all 
restricted to the remaining vertices of $X^i$.)
These properties hold for the following reasons. (C1) holds as every 
vertex deleted from an $X_i$ has been added to $L_e$, whilst (C2) is clear as whenever we deleted vertices we simply 
restricted $G$ and $M$ to the remaining vertices. For (C3), recall that $G_i' \sm M_i'$ was $(c,\eps'')$-robustly $2^k$-universal. Moreover, for all $i\in [t']$ and all $j \in [k]$  
we have $|X^i_j| \geq |X_{i,j}'|/2 \ge c|X_{i,j}'|$, since we ensured that we only deleted
$\theta ''n_1$ vertices from any single $AX'_{i}$ (and at most two from $BX'_i$). Furthermore by 
(B1) we know that $|G^i(v)_=| \ge |(G_i'(v)[ BX'_i])_= | \ge c|G'_i(v)_=|$
for any $v \in X^i$. 
(Also, even if we had arbitrarily deleted a further $\eps n_1$ vertices from $X'_i$ 
when obtaining $X^i$, $G^i$ and $M^i$, these bounds would still hold.)
So $G^i \sm M^i$ satisfies (i) in the definition of a robustly universal complex
(where $X_j^i$ plays the role of $V_j$).
The sets $Z_s$ satisfy (iii) in the definition and so we can find
the required copy of $L$ (even after the deletion of up to $\eps  n_1$ more vertices of $X^i$).
(C4) follows from (\ref{desityHG}) and the fact that $X^i$ was
formed by deleting at most $(\theta'+\theta'') n_1 \ll c' |X_i|$ vertices from $X_i$. Similarly, for (C5) note that
(even after up to $\eps n_1$ more deletions) we have deleted at most $2\theta'' n_1$ vertices from each $X_i$ 
since we split the clusters to form the $X_i$. 
So by~(\ref{splitsizes1}), after these deletions we must have  
\begin{itemize}
\item $\frac{n_1}{2k-3} - 2\theta'' n_1 \leq |X^i_1| \leq \frac{n_1}{2k-3}$,
\item $\frac{n_1}{2k-3} - 2\theta'' n_1 \leq |X^i_2| \leq  \frac{(1-\eta)2n_1}{2k-3}$, and
\item $ \frac{(1-\eta)2n_1}{2k-3} - 2\theta'' n_1 \leq |X^i_j| \leq  \frac{(1-\eta)2n_1}{2k-3}$ for $3 \le j \le k$, 
\end{itemize}
and by~(\ref{splitsizes2}) we must have
\begin{itemize}
\item $
(1-\delta')n_1 - 2(k-1)\theta'' n_1 \le \sum_{j=2}^k  |X_{i,j}| \le (1-\delta)n_1.$
\end{itemize}
Since~$\theta'' \ll \delta \ll \delta' \ll \delta'' \ll \eta$, we deduce that
\begin{itemize}
\item $n_i' \ge \frac{1}{k-1} \left( n_1 \left( 1-\delta' + \frac{1}{2k-3}-2k\theta'' \right)-1 \right)
\ge \frac{(1-\eta)2n_1}{2k-3}$, and 
\item $n'_i \leq \frac{n_1}{k-1} \left( 1-\delta + \frac{1}{2k-3} \right)  \le \frac{(2-\delta)n_1}{2k-3}.$
\end{itemize}
So property (C5) follows.%
    \COMMENT{need~(\ref{splitsizes2}) for upper bound on $n'$. (\ref{splitsizes1}) is not enough}
Finally, (C6) follows from the final step in the construction of $L_e$, in which we added an extra edge to each end of $L_e$ so that (C6) would be satisfied.

\subsection{The supplementary graph} \label{supp}
Roughly speaking, our aim is to find a spanning loose path in each $G^i\setminus M^i$ (and thus in~$H^i$) such that
all these paths together with~$L_e$ form a loose Hamilton cycle in~$H$. So we have to ensure that the complete
$k$-partite $k$-graph on $X^i$ contains a spanning loose path (for this, we will need
$|X^i| \equiv 1\mod k-1$) and we need to join up all the loose paths we find in the~$H^i$.
The purpose of this section is to find the `connecting loose paths' which join up the~$X^i$
in such a way that the divisibility problems are dealt with as well. To do this, we first define a
supplementary hypergraph~$R^*$ whose vertices correspond to the $X^i$. We will show that~$R^*$
is connected and that `along' edges of~$R^*$ we can find our loose paths in~$H$ which join
up all the~$X^i$.

The vertex set of the \emph{supplementary hypergraph~$R^*$} is~$[t']$.
A subset $e\subseteq [t']$ of size at least~2 is an edge of~$R^*$ if there exists an edge
$S_e\in R$ such that for all $j\in S_e$ there are $i_j\in e$ and $\ell_j\in [k]$ with
$X^{i_j}_{\ell_j}\subseteq V_j$ and $e=\{i_j: j\in S_e\}$.
(We fix one such edge~$S_e$ for every $e\in R^*$.)
Then every edge of~$R^*$ has size at most~$k$. We say that $X^i$ \emph{belongs} to an edge $e\in R^*$
if $i\in e$. Similarly, $X^i$ \emph{belongs to} some subhypergraph $R'\subseteq R^*$ if $i\in V(R')$.%
\COMMENT{Had to change things since earlier proof assumed $R^*$ to be uniform which may not be the case}

\begin{lemma} \label{suppgraphconnected}
The supplementary graph~$R^*$ is connected.
\end{lemma}
\nib{Proof.}
Recall that we chose the copies~$A_\ell$ of $\mc{A}_k$ in such a way that the sub-$k$-graph~$A$ of~$R$ induced
by $\bigcup_{\ell=1}^t A_\ell$ is connected. Suppose that $R^*$ is not connected. Let $R^*_1$ be a
component of $R^*$ and let $R^*_2 = R^* -R^*_1$.  Let
$R_1  = \{j\in[m]: X^i_s \subseteq V_j$ for some $i \in V(R^*_1), s \in [k]\}$. So $R_1$ corresponds to the
set of all those clusters which meet some~$X^i$ belonging to~$R^*_1$. Define $R_2$ similarly.%
     \COMMENT{$R_1$ and $R_2$ might intersect, but that does not seem to be a problem}
Then $R_1 \cup R_2 = V(A)$ and thus~$A$ contains some edge~$S$ intersecting both~$R_1$ and~$R_2$.
But then~$S$ corresponds to an edge of~$R^*$ intersecting%
    \COMMENT{Indeed, for each $j\in S$ choose an $i_j\in [t']$ such that $X^{i_j}_\ell\subseteq V_j$
for some~$\ell$. We can choose the $i_j$ such that at least one of them lies in~$V(R^*_1)$ and
at least one lies in~$V(R^*_2)$. (But not all of the $i_j$ need to be distinct.) Then
$\{i_j:j\in S\}$ is an edge of~$R^*$.}
both~$V(R^*_1)$ and~$V(R^*_2)$, a contradiction.
\qed

\medskip

The next lemma shows that within the $X^i$ belonging to an edge of $R^*$, we can find a reasonably
short loose path in $H$ and we may choose (modulo $k-1$) how many vertices this path uses from each $X^i$.
Using the connectedness of~$R^*$, this will allow us to find the connecting loose paths which join up the $X^i$ whilst
having control over the divisibility properties. We shall also insist that the path in
Lemma~\ref{interpath} avoids a number of `forbidden vertices', to enable us to ensure that our connecting loose paths are disjoint, and that the endvertices
of these paths lie in many edges of the relevant $H^i$.

\begin{lemma} \label{interpath}
Suppose that $e\in R^*$ and that for every $i\in e$ there is an integer $t_i$ such that
$0 \leq t_i \leq k-1$ and $\sum_{i\in e} t_i \equiv 1 \mod k-1$. Let $i',i''\in e$ be distinct.
Moreover, suppose that $Z$ is a set of
at most $100(t')^2k^3$ `forbidden' vertices of~$H$. Then in the sub-$k$-graph of $H$
induced by $\bigcup_{i\in e} X^i$ we can find a loose path $L$ with the following properties.%
\COMMENT{previous $P$ was a notation clash}
\begin{itemize}
    \item $L$ contains at most $4k^3$ vertices.
    \item $L$ has an initial vertex $u$ in $X^{i'}$ and a final vertex $v$ in $X^{i''}$.
    \item $|V(L) \cap X^{i}| \equiv t_i \mod k-1$ for each $i\in e$.
    \item $L$ contains no forbidden vertices, i.e. $V(L) \cap Z = \emptyset$.
    \item $u$ lies in at least $|H^{i'}|/(2|X^{i'}|)$ edges of $H^{i'}$, and $v$ lies in at
least $|H^{i''}|/(2|X^{i''}|)$ edges of $H^{i''}$.
\end{itemize}
\end{lemma}
\nib{Proof.}
Recall that in Section~\ref{regsection} we assigned a $k$-partite $k$-complex $G^S$ to every edge $S\in R$
such that (A1)--(A5) are satisfied. To simplify notation, we write~$S$ for the edge $S_e\in R$ corresponding
to~$e$ and suppose that $S=[k]$.
For each $j\in S=[k]$ choose $i_j\in e$ and $\ell_j\in [k]$ such that $X^{i_j}_{\ell_j}\subseteq V_j$
and such that $e=\{i_j: j\in S=[k]\}$. 
To simplify notation we write $Y_j$ for~$X^{i_j}_{\ell_j}\sm Z$, $Y=\bigcup_{j\in [k]}Y_j$
and assume that $i'=i_1$ and $i''=i_k$.
For each $i\in e$ let $J_i$ be the set of all $j\in S=[k]$ with $i_j=i$. So the sets $J_i$ are
disjoint and their union is $[k]$.
Pick some $j\in J_i$ and let $t'_j=t_i$ and $t'_s=0$ for all $s\in J_i\sm \{j\}$.
Our path~$L$ will consist of $t'_j$ vertices from each~$Y_j$ (modulo $k-1$) and thus of
$t_i$ vertices from each~$X^i$ (modulo $k-1$).

Since $G^S$ satisfies (A1)--(A3) and~(A5), Lemma~\ref{restrict} implies that the restriction $G^S[Y]$ is
$\eps'$-regular, with absolute density at least $d(G^S)/2\ge d_a$,
relative density at index $[k]$ at least $d$ and $(G^S)_{\{j\}}[Y]=Y_j$. Furthermore, (A4)
together with the fact that $d(G^S[Y])\ge d(G^S)/2$ imply that 
$$
|M^S[Y]| < |M^S|< \frac{2\nu^{1/4}|G^S|}{c''}\le \theta |G^S[Y]|.$$
Thus Theorem~\ref{robust-universal} implies that we can delete $\theta'|Y_j|$ vertices from each~$Y_j$
to obtain subsets~$Y'_j$ such that $G^S[Y'] \sm M^S[Y']$ is  $(c,\eps'')$-robustly $2^k$-universal, where $Y'=\bigcup_{j\in [k]} Y'_j$.

Now, let $v_j = (k+2)(k-1)+t'_j$. Then $\sum v_j \equiv 1 \mod k-1$ and
so $n' = ((\sum v_i) - 1)/(k-1)$ is an integer. Furthermore, $k(k+2) \leq n' \leq k(k+3)$, and so
$n'/2+1 \leq v_j \leq n'$ for each $j$. Thus by Lemma~\ref{loosepath} we can find a loose path
in the complete $k$-partite $k$-graph on the vertex set~$Y'$, beginning in $Y_1'$, finishing in $Y_k'$
and using $v_j$ vertices from each $Y_j'$. Since $G^S[Y'] \sm M^S[Y']$ is $(c, \eps'')$-robustly $2^k$-universal,
we can find such a loose path $L$ in $G^S[Y']\sm M$ and hence in $H-Z$.
(Indeed, we can do this by finding the complex $L^\leq$, which has maximum vertex degree at most $2^k$.
Note that we use the definition with $J=G'$ in (i)).
Note that~$L$ contains at most $k(k-1)(k+3) \leq 4k^3$ vertices. 

To see that we can insist on the final condition of the lemma, recall that $d(H^i)\ge c'$ by~(C4).
Thus for all $j\in [k]$ at least $c'|X^i_j|/2$ vertices of $X^i_j$ lie in at
least $|H^i|/(2|X^i_j|)$ edges of $H^i$, and so we may restrict the initial and final vertices of $L$
to these sets of vertices (minus the vertices of $Z$) by~(iii) in the definition of
robust universality. \qed

\subsection{Constructing the loose Hamilton cycle}\label{linking}
As discussed before, our Hamilton cycle in~$H$ will consist of~$L_e$ and paths in each~$H^i$
as well as paths connecting the~$X^i$. However, we need to make sure that all these paths join up
nicely, motivating the following definition. Suppose $L$ is a path in some $k$-graph $K$ 
with initial vertex $x'$ and final vertex $y'$.
Also, let $I, F\subseteq V(K)\sm V(L)$ be disjoint sets of size $k-2$. Then $L^*=I \cup F \cup V(L)$ is a
\emph{prepath}. Note that $L^*$ is not (the vertex set of) a $k$-graph, but that if we
can find vertices $x, y \in V(K)\sm L^*$ such that $\{x,x'\} \cup I, \{y, y'\} \cup F \in K$,
then adding $x$ and $y$ to $L^*$ gives another path. We refer to all such vertices $x \in V(K)$
as \emph{possible initial vertices} of~$L^*$ and to all such vertices $y \in V(K)$ as \emph{possible final vertices}. 
If $L$, $L'$ and $L''$ are disjoint loose paths, $I,F,x,y$ are as before,
$x$ is also the final vertex of~$L'$ and $y$ is also the initial vertex of $L''$ then $I$ and~$F$ together with
$L',L,L''$ form a single loose path, illustrating how we shall join paths together. 

We start by converting our exceptional path $L_e$ into a prepath. Recall that $|V(L_e)|<\eta n$ and
that the initial vertex $x_e$ of~$L_e$ and its final vertex $y_e$ satisfy~(C6). Let $a\in [t']$
and $u_a\in [k]$ be such that $x_e\in X^{a}_{u_a}$. Pick any $u'_a\in [k]$ with $u_a\neq u'_a$.
(C4) and (C6) together imply that there is a set
$I_0\subseteq X^{a}\setminus (X^{a}_{u_a}\cup X^{a}_{u'_a})$
for which $X^{a}_{u'_a}$ contains at least $c|X^{a}|$ vertices~$v$
which form an edge of $H^a$ together with $I_0\cup \{x_e\}$. Let $I'_0\subseteq X^{a}_{u'_a}$ be such
a set of vertices. Similarly, letting $b\in [t']$, $u_b\neq u'_b\in [k]$ be such that
$y_e\in X^b_{u_b}$, there is a set $F_0\subseteq X^{b}\setminus (X^{b}_{u_b}\cup X^{b}_{u'_b}\cup I_0)$
for which $X^{b}_{u'_b}$ contains at least $c|X^{b}|$ vertices~$v$
which form an edge of $H^b$ together with $F_0\cup \{y_e\}$. Let $F'_0\subseteq X^{b}_{u'_b}$ be such
a set of vertices. Let $L^*_e$ be the prepath $I_0\cup F_0\cup V(L_e)$. Then $I'_0$ is a set of possible
initial vertices of $L_e^*$ and $F'_0$ is a set of possible final vertices. (We do not remove $I_0$ from
$X^a$ and $F_0$ from $X^b$ at this stage.)

Since by Lemma~\ref{suppgraphconnected} the supplementary graph $R^*$ is connected, we can find a walk~$W$
from $b$ to $a$ in $R^*$ such that every $i\in [t']=V(R^*)$ appears as an initial, link or final
vertex in~$W$ (these vertices were defined in Section~\ref{walks}) and such that $W$ has length $\ell \leq (t')^2$.
Let $e_1,\dots,e_\ell$ be the edges of this walk, let $r_1 = b, r_{\ell+1} = a$, and
let $r_2, \dots, r_{\ell}$ be the link vertices of the walk. For each $i\in [t']$, let
$d_i = |\{j\in [\ell+1]: r_j = i\}|$, that is, the number of times~$i$ appears as an initial, link or final vertex
in $W$. So $d_i > 0$ for every $i$ and $\sum d_i = \ell+1 $.

Our next aim is to apply Lemma~\ref{interpath}
to each edge $e_j$ in order to find a loose path $L_j$ in $H$, which we will extend to a prepath $L_j^*$
with many possible initial vertices in $X^{r_j}$ and many possible final vertices in $X^{r_{j+1}}$.
We shall do this for each $e_1,\dots,e_\ell$ in turn. So suppose that $s\in [\ell]$ and that
for all $j=1,\dots,s-1$ we have defined loose paths~$L_j$ in~$H$ as well as sets $I_j,F_j$ extending~$L_j$
to a prepath~$L^*_j$ which satisfy the following properties:
\begin{itemize}
\item[(D1)] $L_j$ lies in the sub-$k$-graph of $H$ induced by $\bigcup_{i\in e_j} X^i$
and contains at most $4k^3$ vertices.
\item[(D2)] The initial vertex~$x_j$ of~$L_j$ lies in $X^{r_j}$ and its final vertex~$y_j$ lies
in $X^{r_{j+1}}$.
\item [(D3)] $I_j\subseteq X^{r_j}$ and $F_j\subseteq X^{r_{j+1}}$.
\item[(D4)] There is a set $I'_j\subseteq X^{r_j}$ of at least $c|X^{r_j}|$
possible initial vertices for $L^*_j$. Similarly, there is a set $F'_j\subseteq X^{r_{j+1}}$
of at least $c|X^{r_{j+1}}|$ possible final vertices for $L^*_j$.
\item[(D5)] All the prepaths $L^*_e,L^*_1,\dots,L^*_{s-1}$ are disjoint.
\item[(D6)] For each $i\in [t']$ and all $j=0,\dots,s-1$ let $X^i(j) = X^i\sm (V(L_1)\cup\dots\cup V(L_j))$,
where $X^i(0)=X^i$. For each $j\in [s-1]$ set $t_i(j) = |X^i(j-1)|+d_i$.
Then for every $i\in e_j$ with $i\neq r_{j+1}$ we have $|V(L_j)\cap X^i|\equiv t_i(j) \mod k-1$.
Moreover $|V(L_j)\cap X^{r_{j+1}}|\equiv 1-\sum_{i\in e_j,\ i\neq r_{j+1}} t_i(j) \mod k-1$.
\end{itemize}
Let us now show how to find $L_s$, $I_s$ and $F_s$. Apply Lemma~\ref{interpath} with $e=e_s$,
$i'=r_s$, $i''=r_{s+1}$ and with $Z=L^*_1\cup\dots L^*_{s-1}\cup I_0\cup F_0$
to find a loose path $L_s$ which satisfies (D1), (D2), (D6) and is disjoint from
$L^*_e,L^*_1,\dots,L^*_{s-1}$. Moreover, the initial vertex $x_s$
of~$L_s$ lies in at least $|H^{r_s}|/(2|X^{r_s}|)$ edges of $H^{r_s}$, and the final vertex
$y_s$ of $L_s$ lies in at least $|H^{r_{s+1}}|/(2|X^{r_{s+1}}|)$ edges of $H^{r_{s+1}}$.
We can now use the latter property to choose sets $I_s$ and $F_s$ which extend $L_s$ to a prepath~$L^*_s$
satisfying~(D3)--(D5). The argument for this is similar to that for the extension of~$L_e$ to $L^*_e$.
Altogether this shows that we can find prepaths $L^*_1,\dots,L^*_\ell$ satisfying (D1)--(D6).

For each $i \in [t']$ we let $j_i$ be the maximal integer such that $i \in e_{j_i}$.
Thus $X^i(\ell) = X^i(j_i) = X^i(j_i-1) \sm V(L_{j_i})$ by~(D1). But if $i \neq r_{\ell+1}$ then (D5) and~(D6)
together imply that
$$|V(L_{j_i}) \cap X^i(j_i-1)|=|V(L_{j_i}) \cap X^i| \equiv t_i(j_i) \equiv |X^i(j_i-1)|+d_i \mod k-1$$
and so $|X^i(\ell)| \equiv -d_i \mod k-1$. 
We claim that this also holds if $i = r_{\ell+1}$. 
To see this, recall that since $L_j$ is loose, we have $|V(L_j)|\equiv 1\mod k-1$ for each
$j\in [\ell]$. Hence
\begin{eqnarray*}
|X^{r_{\ell+1}}(\ell)| & = & |V(H) \sm V(L_e)|-\sum_{j\in [\ell]} |V(L_j)|-
\sum_{i\in [t'],\ i\neq r_{\ell+1}}|X^i(\ell)|\\
& \stackrel{(\ref{exceppathsize})}{\equiv} & -1-\ell+\sum_{i\in [t'] \ i\neq r_{\ell+1}} d_i 
\equiv -d_{r_{\ell+1}} \mod k-1
\end{eqnarray*}
as $\ell+1=\sum_{i\in [t']} d_i$. Let $Y^i=X^i\sm (L^*_e\cup L^*_1\cup\dots\cup L^*_\ell)$.
Since by~(D3) for each $i\in [t']$ there are exactly $2(k-2)d_i$
vertices of $X^i$ which lie in $L^*_e,L^*_1,\dots,L^*_\ell$ but not in $L_e,L_1,\dots,L_\ell$,
this in turn implies that
\begin{equation}\label{Yi}
|Y^i|\equiv -d_i-2(k-2)d_i \equiv d_i \mod k-1.
\end{equation}
Let $x_{\ell+1}=x_e$, $y_0=y_e$, $L^*_0=L^*_e$, $I_{\ell+1}=I_0$ and $I'_{\ell+1}=I'_0$.
In order to complete our prepaths $L^*_0,\dots,L^*_\ell$ to a Hamilton cycle we wish
to choose~$d_i$ disjoint loose paths $L^i_1,\dots,L^i_{d_i}$ within each $H[Y^i]$ 
which together contain all the vertices in $Y^i$ and which `connect' successive prepaths $L_j^*$.
We achieve this as follows.
Let $J_i$ be the set of all $j\in [\ell+1]$ with $r_j=i$. So $J_i$ is the set of positions
at which~$i$ occurs as an initial, final or link vertex in our walk~$W$ and $|J_i|=d_i$.
Let $j_1\le \dots\le j_{d_i}$ be the elements of~$J_i$. 
Then we choose the $L^i_s$ ($s\in [d_i]$)
in such a way that the initial vertex of $L^i_s$ lies in $F'_{j_s-1}$ and its final vertex
lies in $I'_{j_s}$, all the $L^i_s$ are disjoint and together they cover all the vertices in $Y^i$.
To see that this can be done, first note that 
$|X^i\sm Y^i|\le \ell (4k^3+2(k-2))+2(k-2) \ll \eps n_1$. So 
using Lemma~\ref{loosepath} together with~(C5) and~(\ref{Yi}) it is easy to check 
that the complete $k$-partite $k$-graph
on~$Y^i$ contains such paths (e.g.~first choose $L^i_1,\dots,L^i_{d_i-1}$, each consisting
of precisely 2 edges, and then apply~(C5) and Lemma~\ref{loosepath} to find a loose path~$L^i_s$
containing all the remaining vertices of~$Y^i$). Now
(C3) and (D4) together imply that 
$G^i[Y^i]\sm M^i[Y^i]$ contains the $k$-complexes induced by these paths
(i.e.~it contains $(L^i_1)^\le,\dots,(L^i_{d_i})^\le$). But this means that we can find the
required paths $L^i_1,\dots,L^i_{d_i}$ in each $H[Y^i]$.

Finally, for each $s\in [d_i]$ write $L'_{j_s}$ for $L^i_s$ and $x'_{j_s}$ for its initial and $y'_{j_s}$ for its
final vertex (where $j_s$ is as defined in the previous paragraph). To obtain our Hamilton cycle in~$H$
we first traverse $L_0=L_e$, then we use the edge $F_0\cup \{y_0,x'_1\}$ in order to move to
the initial vertex $x'_1$ of $L'_1$. (This is possible since $x'_1\in F'_0$.) Now we traverse
$L'_1$ and use the edge $I_1\cup \{y'_1,x_1\}$ to get to~$x_1$. (Again, this is possible
since $y'_1\in I'_1$.) Next we traverse $L_1$ and use the edge $F_1\cup \{y_1,x'_2\}$
to move to~$x'_2$. We continue in this way until we have reached the initial vertex $x_{\ell+1}=x_e$ of
$L_0=L_e$ again. (So in the last step we traversed $L'_{\ell+1}$ and used the edge
$I_{\ell+1}\cup \{y'_{\ell+1},x_{\ell+1}\}$.) This completes the proof of Theorem~\ref{main}. \endproof

\medskip

\noindent
{\footnotesize
Peter Keevash, Richard Mycroft, School of Mathematical Sciences, Queen Mary, University of London, Mile End Road,
London, E1 4NS, United Kingdom, 
\{{\tt p.keevash,r.mycroft}\}{\tt @qmul.ac.uk }

\smallskip

\noindent
Daniela K\"uhn, Deryk Osthus, School of Mathematics,
University of Birmingham, Birmingham, B15 2TT, United Kingdom, 
 \{{\tt kuehn,osthus}\}{\tt @maths.bham.ac.uk } }








\end{document}